# Stochastic volatility model with long memory for water quantity-quality dynamics

Hidekazu Yoshioka[a,*] and Yumi Yoshioka[b]

[a]Japan Advanced Institute of Science and Technology, 1-1 Asahidai, Nomi 923-1292, Japan
[b]Gifu University, 1-1 Yanagido, Gifu 501-1193, Japan
[*]Corresponding author: yoshih@jaist.ac.jp, ORCID: 0000-0002-5293-3246

**Abstract:** Water quantity and quality are vital indices for assessing fluvial environments. These indices are highly variable over time and include sub-exponential memory, where the influences of past events persist over long durations. Moreover, water quantity and quality are interdependent, with the former affecting the latter. However, this relationship has not been thoroughly studied from the perspective of long-memory processes, which this paper aims to address. We propose applying a new stochastic volatility model, a system of infinite-dimensional stochastic differential equations, to describe dynamic asset prices in finance and economics. Although the stochastic volatility model was originally developed for phenomena unrelated to the water environment, its mathematical universality allows for an interdisciplinary reinterpretation: river discharge is analogous to volatility, and water quality to asset prices. Moreover, the model's infinite-dimensional nature enables the analytical description of sub-exponential memory. The moments and autocorrelations of the model are then obtained analytically. We mathematically analyze the stochastic volatility model and investigate its applicability to the dynamics of water quantity and quality. Finally, we apply the model to real time-series data from a river in Japan, demonstrating that it effectively captures both the memory and the correlation of water quality indices to river discharge. This approach, grounded in infinite-dimensional stochastic differential equations, represents a novel contribution to the modelling and analysis of environmental systems where long memory processes play a role.

**Fundings** This study was supported by the Japan Society for the Promotion of Science (KAKENHI No. 22K14441) and the Japan Science and Technology Agency (PRESTO No. JPMJPR24KE).

**Competing interests** The authors have no relevant financial or non-financial interests to disclose.

**Acknowledgments** The authors would like to express their gratitude to Dr. Ikuo Takeda of Shimane University for providing a part of the water quality data used in this paper.

**Declaration of generative AI in scientific writing** The authors did not use generative AI for scientific writing of this manuscript.

1. **Introduction**

**1.1 Study background**

Assessing water quantity and quality in rivers is a core issue for fostering sustainable and harmonious coexistence between nature and humans. Examples of such efforts globally include, but are not limited to, hydrological connectivity in river systems (Durighetto et al., 2023; Zanetti et al., 2024)[1,2], surveys of international rivers under severe anthropogenic pressure (Ding et al., 2024)[3], flood hazard management in the context of climate change (Bahramloo et al., 2024)[4], investigations into potential threats to biodiversity loss in river environments (Rai et al., 2024)[5], fish habitat assessments in riverscapes (Hansen et al., 2024)[6], and water quality analyses in lakes (Gao et al., 2024)[7].

The physical, chemical, and biological dynamics of river environments are highly variable and subject to random fluctuations. Consequently, it is reasonable to conceptualize these dynamics as stochastic processes, characterized by continuous and discrete time-series data that fluctuate randomly, such as environmental variables including rainfall intensity and river discharge (Doukhan, 2018; Volpi et al., 2024)[8,9]. Particularly, stochastic differential equations (SDEs) (Øksendal and Sulem, 2019; Pascucci, 2024)[10,11] are well-suited for balancing mechanistic and phenomenological modeling. SDEs play a central role in ecological modeling, where randomness significantly affects system stability (Mandal et al., 2024; Paul et al., 2025)[12,13]. They have also been effectively used to investigate disease propagation (Mohammad et al., 2024; Đorđević and Dahl, 2024)[14,15] and asset price dynamics and green financing in economics (Hambel et al., 2024; Zhang and Ren, 2024)[16,17]. Coastal sand dynamics (Ramakrishnan et al., 2024; Vinent et al., 2021)[18,19] and salt intrusion in tidal marshes (Dijkstra et al., 2023)[20] are additional examples of SDE applications.

The rate at which autocorrelation decays characterizes the forgetting timescale of past events. Traditionally, this decay has been assumed to follow an exponential pattern, or Markovian dynamics, due to its analytical tractability in applications (Calvani and Perona, 2023; Latella et al., 2024; Olson et al., 2021; Yoshioka and Yoshioka, 2023a; Yoshioka and Yoshioka, 2023b)[21-25]. However, empirical hydrological data, including discharge, water quality, and hydrological signatures, often exhibit sub-exponential decay, indicative of longer memory and non-Markovian behavior (Fatni et al., 2024; Maftei et al., 2016; Guo et al., 2024a; Pizarro et al., 2024; Rahmani and Fattahi, 2024)[26-30]. Therefore, employing non-Markovian models with long memory is a reasonable strategy for addressing issues related to water quantity and quality.

Some non-Markovian processes can alternatively be expressed as a superposition (i.e., suitable integration) of Markovian processes, albeit at the cost of increased dimensionality, which may extend from finite to infinite dimensions (e.g., Abi Jaber, 2019; Cuchiero and Teichmann, 2020)[31,32]. This mathematical technique is called Markovian lifts or Markovian embedding, depending on the research domain, with slight differences in terminology. Markovian lifts have been employed in various fields, including physics (Kanazawa and Sornette, 2024; Wiśniewski and Spiechowicz, 2024)[33,34] and economics (Bondi et al., 2024; Damian and Frey, 2024; Dupret and Hainaut, 2024)[35-37]. The Markov chain approximation of the Volterra process is an example of Markovian lifts (Yang et al., 2024)[38]. In computational fluid dynamics, the principle of decomposing long memory into a sum of shorter memories has been recently explored (Guo et al., 2024b)[39], suggesting the applicability of Markovian lifts across diverse research areas. Moreover, the use of numerous hypothetical reservoir cells has been shown to be effective in physics-informed machine learning models (He et al., 2024)[40].

Despite their mathematical complexity, Markovian lifts have practical implications in the context of autocorrelation. They decompose sub-exponential functions into integrals of exponential functions, with exponents distributed according to a probability measure (e.g., Fasen and Klüppelberg, 2005)[41]. This property has been exploited in modeling long-memory dynamics in water quantity and quality (Yoshioka and Yoshioka, 2024a; Yoshioka and Yoshioka, 2024b)[42,43]. However, these studies have focused on either water quantity or water quality in isolation. The coupled dynamics of water quantity and quality which interact in real river environments, remain unexplored from the perspective of long-memory modeling via Markovian lifts. Addressing this gap requires innovative approaches to integrate the two different dynamics while preserving the analytical tractability of Markovian lifts and accurately capturing their statistical characteristics.

Stochastic volatility models in asset dynamics have demonstrated that long-memory processes, based on a superposition approach, can be coupled with another stochastic differential equation (SDE) without compromising overall tractability (Barndorff-Nielsen and Stelzer, 2013)[44]. In such models, volatility is represented as the fluctuating diffusion coefficient of an SDE. The key to this approach lies in its affine nature (e.g., Duffie et al., 2003)[45], which allows statistics such as mean, variance, and autocovariance to be obtained in closed form. Stochastic volatility models with long-memory have been extensively studied, focusing on the roughness of their sample paths (Bayer et al., 2021; de Truchis et al., 2024; Wu et al., 2022)[46-48]. While such models have been applied to river discharge modeling (Wang et al., 2023)[49], they have yet to incorporate water quality dynamics.

Water quality models for river environments based on SDEs have addressed various aspects, including the influence of suspended sediment concentrations on fluvial processes (Jing et al., 2020)[50], dissolved oxygen and biological oxygen demand (Mansour, 2023)[51], antibiotic bacterial dynamics (Gothwal and Thatikonda, 2020)[52], nutrient loads with a focus on regime shifts in river environments (Park and Rao, 2014)[53], and phycocyanin as a pigment indicator of cyanobacteria (Carpenter and Brock, 2024)[54]. However, most existing SDE models for water quality dynamics are Markovian, and therefore do not adequately account for memory effects, which are known to influence a wide variety of water quality indices (WQIs) (Chong et al., 2023; Rahmani and Fattahi, 2024; Spezia et al., 2021; Yoshioka and Yoshioka, 2024b)[55-57,43]. Another critical issue is the concentration-discharge relationship, which expresses the concentration of WQIs as a function of discharge. This relationship has been studied for various indices (Qin et al., 2024; Tunqui Neira et al., 2021; Wichman et al., 2024; Wymore et al., 2023; Zhan et al., 2022)[58-62] and quantifies the correlation between water quality and quantity. Despite its importance, this relationship has not been adequately incorporated into existing SDE models, leaving a gap in the development of realistic stochastic models for coupled water quantity-quality dynamics.

These challenges motivated this study, which seeks to address the aforementioned gaps by conceptualizing water quantity and quality as analogous to asset prices and volatility, respectively. The aims and contributions of this study are detailed below.

**1.2 Aim and contribution**

The aims of this study are twofold: (i) the proposal and analysis of a new coupled water quantity-quality model within the framework of affine processes, and (ii) the application of the model to real-world data. The study's contributions to achieving these aims are detailed below.

**(i)   Modeling and analysis**

The proposed model is formulated as two coupled superposed processes representing river discharge and a WQI. It is based on the superposition approach introduced by Barndorff-Nielsen (2001)[63], in which a long-memory process is generated by superposing infinitely many "small", measure-valued Ornstein–Uhlenbeck (supOU) processes with exponential memory decay. These small processes, characterized by different reversion timescales, are combined via stochastic integration using appropriate random measures. Initially, small processes are chosen as linear SDEs but are later generalized to affine processes, with nonlinear extensions being examined (Iglói and Terdik, 2003)[64]. The supOU process is particularly tractable, making it useful in applications such as option pricing models with memory (Leonenko et al., 2024)[65] and in exploring statistical properties such as large deviations and fractal properties (Grahovac et al., 2019; Grahovac and Kevei, 2025)[66,67]. The success of the superposition approach lies in the affine form of the coefficients of the infinitesimal generators associated with these small processes (e.g., Chevalier et al., 2022; Motte and Hainaut, 2024; Yoshioka et al., 2023)[68-70], which enable the derivation of the generalized Riccati equation governing the moment-generating function of the superposed process.

In this paper, we realize a stochastic volatility model for coupled water quantity and quality dynamics by first formulating a jump-driven supOU process to represent river discharge. The temporal evolution of the WQI is then modeled as the product of a deterministic seasonal function and the exponential of an SDE, with drift and diffusion modulated by the supOU process. This framework allows for the representation of both positive and negative correlations between river discharge and the WQI (**Figure 1**). We refer to this coupled model as the superposed Ornstein–Uhlenbeck stochastic volatility (supOUSV) model. The supOUSV model naturally incorporates the concentration-discharge relationship without compromising analytical tractability. As shown in **Figure 1**, the model of discharge and WQI is coupled in a unidirectional manner so that the dynamics of WQI are influenced by those of discharge assuming that the quality of water does not affect its quantity. The concentration-discharge relationship in the proposed model essentially results from this assumption. We derive a closed-form expression for the autocorrelation of the WQI, along with its tail estimates. Moreover, the approach facilitates the derivation of a generalized Riccati equation for the moment-generating function, ensuring boundedness of the solution. Generalized Riccati equations are pivotal for analyzing infinite-dimensional stochastic volatility models (Cox et al., 2022; Friesen and Karbach, 2024)[71,72]. The specific mathematical structure of the supOUSV model allows for analytically obtaining several key statistics, which is advantageous in both theory and applications, as demonstrated in a case study explained below.

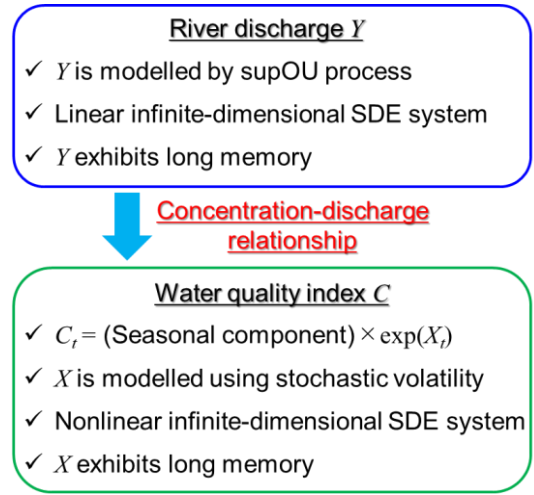

**Figure 1.** Conceptual diagram of the proposed model.

**(ii) Model application**

We apply the proposed supOUSV model to real river discharge and WQI data at a study site in the midstream section of the Hii River, a Class A river in Japan. Hourly river flow time-series data are publicly available for the site, and almost weekly sampled data of multiple WQIs have been collected over 30 years. The empirical data indicate that both discharge and logarithmic WQIs include long-memory characteristics. The model is fitted to the data based on least-squares and moment-matching methods, enabling computational investigation of the proposed model. The results show that the model effectively captures both the long-memory dynamics and the discharge dependence of WQIs. Additionally, we analyze the concentration-discharge relationship, focusing on the hysteresis patterns of the WQI (Mazilamani et al., 2024)[73], as represented by the model. In summary, this study contributes to the formulation, analysis, and application of a novel mathematical model for long-memory processes, advancing the understanding of coupled water quantity and quality dynamics.

## 1.3 Organization of this paper

The remainder of this paper is organized as follows. **Section 2** describes our model. **Section 3** mathematically analyzes the model by focusing on statistics. **Section 4** describes the application of the model to river discharge and WQIs at a study site in Japan. **Section 5** concludes this paper and presents its perspectives. **Appendices** contain results concerning the finite-dimensional version of the proposed model (**Section A. 1**), proofs of propositions (**Section A. 2**), and auxiliary results (**Sections A. 3–A. 4**).

## 2. Stochastic volatility model
## 2.1 Overview

We work on a complete probability space $(\Omega, \mathbb{F}, \mathbb{P})$ as assumed in stochastic analysis (e.g., Definition 1.1.7 in Pascucci (2024))[11] ($\Omega$ is collection of events, $\mathbb{F}$ is filtration, and $\mathbb{P}$ is probability). The primary variables in the model are river discharge and WQI. River discharge is represented as a supOU process driven by a pure-jump subordinator as in previous studies (e.g., Yoshioka and Yoshioka, 2024a)[42] (**Figure 1**). The WQI, the focus of this study, is modeled as a positive scalar process expressed as the product of a seasonal part and an exponential superposition of the stochastic volatility SDE.

**2.2 River discharge**

River discharge (referred to as discharge in the sequel) in this study represents the volume of water passing through a cross-section of a river per unit of time. The discharge at a station located within the river is modeled as a nonnegative, continuous-time scalar stochastic process, $Y = (Y_t)_{t \in \mathbb{R}}$, where $t \in \mathbb{R}$ is time. Assuming a stationary state (e.g., Botter, 2010)[74], the discharge is modeled as a supOU process (Barndorff-Nielsen, 2001; Barndorff-Nielsen and Stelzer, 2013)[44,63]:

$$Y_t = \int_{s=-\infty}^{s=t} \int_{z=0}^{z=+\infty} \int_{r=0}^{r=+\infty} z e^{-r(t-s)} N(ds, dz, dr), \ t \in \mathbb{R}. \tag{1}$$

Here, $z > 0$ is the jump size at each jump (i.e., flood event), $r > 0$ is recession rate of a flood event and is distributed according to the probability measure $\pi = \pi(dr)$, and $N$ is a Poisson random measure on $\mathbb{R} \times (0, +\infty) \times (0, +\infty)$ with the compensator $\pi(dr)\nu(dz)ds$, where $\nu$ is the Lévy measure of a pure-jump subordinator (i.e., a pure-jump Lévy process with positive jumps). The two measures $\nu$ and $\pi$ control the size-frequency relationship and the recession of jumps, respectively. Here, $N$ is extended to $\mathbb{R}$ as in Barndorff-Nielsen and Stelzer (2013)[44] to deal with the stationary process.

Physically, the supOU process represents a situation where each flood event, assumed to occur randomly in time, has its own recession rate. The heterogeneity in recession rates originates from the interplay between antecedent water storage in the catchment and rainfall intensity (e.g., Hameed et al., 2023; Lee et al., 2023)[75,76]. This heterogeneity was previously described as a finite sum of independent Ornstein–Uhlenbeck processes by Botter (2010)[74]. The supOU process generalizes this concept to infinite dimensions. The long memory of discharge arises from catchment heterogeneity, where both fast (surface water) and slow runoffs (groundwater) coexist (Mudelsee, 2007; Di Dato et al., 2023)[77,78].

For technical reasons, we assume the following conditions for $\nu$ and $\pi$:

$$\int_0^{+\infty} z^k \nu(dz) < +\infty \ (k = 1, 2, 3, ...) \tag{2}$$

and

$$\int_0^{+\infty} r^{-1} \pi(dr) < +\infty. \tag{3}$$

Assumption (2) means that the jumps in the supOU process admit moments of arbitrary order. This condition is not restrictive for applications to discharge, e.g., a tempered stable model: $\nu(dz) = a_1 e^{-a_2 z} z^{-(1+a_3)} dz$ with $a_1, a_2 > 0$ and $a_3 < 1$ (e.g., Yoshioka and Yoshioka, 2024a)[42]. The tempered stable model includes the finite- ($a_3 < 0$, $\int_0^{+\infty} \nu(\mathrm{d}z) < +\infty$) and infinite-activity case ($a_3 \in [0,1)$, $\int_0^{+\infty} \nu(\mathrm{d}z) = +\infty$) in a unified framework. The assumption (3) ensures that the distribution of the recession rate is not singular at the origin $r = 0$. This condition is satisfied for the gamma distribution $\pi(\mathrm{d}r) = \gamma_{\alpha,\beta}(r)\mathrm{d}r = \frac{1}{\Gamma(\alpha)\beta^\alpha} r^{\alpha-1} e^{-r/\beta} \mathrm{d}r$ with $\alpha > 1$ and $\beta > 0$.

### 2.3 Water quality index
#### 2.3.1 Seasonal decomposition

We propose an SDE model for the WQI of a river. The WQI $C = (C_t)_{t \in \mathbb{R}}$, representing the concentration of a solute, is assumed to be nonnegative and consists of seasonal and stationary parts (Yoshioka and Yoshioka, 2024b)[43]:

$$C_t = \bar{C} e^{S_t + X_t}, \quad t \in \mathbb{R}. \tag{4}$$

Here, $\bar{C} > 0$ is a constant; $S = (S_t)_{t \in \mathbb{R}}$ is a deterministic, periodic, real-valued coefficient; and $X = (X_t)_{t \in \mathbb{R}}$ is the logarithm of the stationary part that encodes the (hidden) long memory of WQI. Some "delay" or "memory" in the river discharge is encoded in (1), which propagates to the WQI through the relationship (4). Under this formulation, water quality investigations aim to identify the two processes $X$ and $Y$, as well as the periodic part $S$.

We assume annual periodicity ($T = 1$ year or $T = 365.25$ days, adjusted for leap years) and expressed the seasonal part as follows:

$$S_t = \sum_{i=1}^{n} A_i \sin\left(\frac{2i \mathrm{PI} t}{T} + B_i\right), \quad t \in \mathbb{R} \tag{5}$$

with some $n \in \mathbb{N}$ and real parameters $A_i$ and $B_i$, and the circle ratio PI=3.1415.... The memory of the WQI is assumed to be encoded in $X$, which is described using a stochastic volatility model. We therefore focus on mathematically modeling $X$. Without loss of generality, we set $\mathbb{E}[X_t] = 0$ for all $t \in \mathbb{R}$, where $\mathbb{E}$ is expectation.

#### 2.3.2 Classical stochastic volatility model

The core of the proposed model is the SDE for the auxiliary scalar process $X$. A classical stochastic volatility model takes the form (e.g., Barndorff-Nielsen and Stelzer, 2013)[44]: for $t \in \mathbb{R}$,

$$\mathrm{d}X_t = -R X_t \mathrm{d}t + \sigma \sqrt{R Y_t} \mathrm{d}B_t \tag{6}$$

or equivalently

$$X_t = \sigma \int_{-\infty}^{t} \sqrt{R Y_s} e^{-R(t-s)} \mathrm{d}B_s, \tag{7}$$

where $B = (B_t)_{t \in \mathbb{R}}$ is a (double-sided) standard Brownian motion independent of the Poisson random measure $N$, $R > 0$ is the reversion rate, and $\sigma > 0$ is the noise intensity. The last term in (6) is defined in the Itô's sense. The SDE in (6) is driven by a noise term with stochastic volatility $\sigma\sqrt{RY_t}$.

A key limitation of the SDE in (6) is that it results in an exponential autocorrelation. Indeed, from (7), for any $h \geq 0$ and $t \in \mathbb{R}$ (noting that $\mathbb{E}[X_t] = 0$):

$$\mathbb{E}[X_t X_{t+h}] = \frac{\sigma^2 \bar{Y}}{2} e^{-Rh}. \tag{8}$$

Here, we used the Itô's isometry (e.g., Theorem 10.1.6 of Pascucci (2024))[11] and the stationarity of $Y$, with the average denoted by $\bar{Y}$. The autocorrelation decays exponentially as shown in (8). Another drawback of (6) is that $X$ does not correlate with $Y$: for $t \in \mathbb{R}$,

$$\mathbb{E}[X_t Y_t] = 0. \tag{9}$$

This lack of correlation restricts the model's applicability, as it assumes no interaction between discharge and water quality.

### 2.3.3 SDE model

The two drawbacks, exponential and hence short memory and no correlation to discharge, are overcome by the proposed supOUSV model. For $t \in \mathbb{R}$,

$$X_t = \underbrace{\int_{R=0}^{R=+\infty} Z_t(\mathrm{d}R)}_{\text{Superposition to generate long memory}} \tag{10}$$

with $Z = \left(Z_t(\mathrm{d}R)\right)_{t \in \mathbb{R}}$ being a measure-valued small process governed by a formal SDE for $R > 0$:

$$\mathrm{d}Z_t(\mathrm{d}R) = \underbrace{-R\left(Z_t(\mathrm{d}R) - \mu(Y_t - \bar{Y})\rho(\mathrm{d}R)\right)\mathrm{d}t}_{\text{Reversion to the mean 0}} + \underbrace{\sigma\sqrt{RY_t}B(\mathrm{d}R, \mathrm{d}t)}_{\text{Volatility-driven fluctuation}}, \tag{11}$$

also expressed as

$$Z_t(\mathrm{d}R) = \int_{s=-\infty}^{s=t} e^{-R(t-s)}\left\{R\mu(Y_s - \bar{Y})\rho(\mathrm{d}R)\mathrm{d}s + \sigma\sqrt{RY_s}B(\mathrm{d}R, \mathrm{d}s)\right\}. \tag{12}$$

Here, $B(\mathrm{d}R, \mathrm{d}t)$ is a space-time Gaussian random measure (e.g., Gomez et al., 2023)[79] with covariance $\mathbb{E}[B(\mathrm{d}R, \mathrm{d}t)B(\mathrm{d}R', \mathrm{d}t')] = \delta(R - R')\delta(t - t')\rho(\mathrm{d}R)\rho(\mathrm{d}R')\mathrm{d}t\mathrm{d}t'$, $\rho(\mathrm{d}R)$ is a probability measure of the distributed reversion $R > 0$, $\mu \in \mathbb{R}$ modulates the dependence of the water quality on discharge, and $\delta$ is the Dirac delta. We assume that the random measures $B$ and $N$ are mutually independent. Assuming a stationary state, the system (10) is rewritten as

$$X_t = \int_{R=0}^{R=+\infty}\int_{s=-\infty}^{s=t} e^{-R(t-s)}\left\{R\mu(Y_s - \bar{Y})\rho(\mathrm{d}R)\mathrm{d}s + \sigma\sqrt{RY_s}B(\mathrm{d}R, \mathrm{d}s)\right\}. \tag{13}$$

The supOUSV model is a system containing (1) and (13), and the WQI is obtained from (4).

## 3. Mathematical analysis
### 3.1 Moments and autocorrelation

### 3.1.1 supOU process

The moments (and hence cumulants) and autocorrelation of the supOUSV model can be obtained explicitly under certain conditions. The cumulants and autocorrelation of the supOU process are as follows (e.g., Barndorff-Nielsen and Stelzer, 2013; Yoshioka, 2022)[44,80]. Here, $\mathbb{V}$, $\mathbb{S}$, and $\mathbb{K}$ are variance, (unnormalized, and hence not divided by using $\mathbb{V}$) skewness and kurtosis, respectively, and $\mathrm{AC}_Y(h)$ is the autocorrelation of $Y$ with the time lag $h \geq 0$:

$$\mathbb{E}[Y_t] = \overline{Y} = M_1 \int_0^{+\infty} \frac{\pi(\mathrm{d}r)}{r}, \quad \mathbb{V}[Y_t] = \overline{V} = \frac{M_2}{2} \int_0^{+\infty} \frac{\pi(\mathrm{d}r)}{r}, \tag{14}$$

$$\mathbb{S}[Y_t] = \frac{M_3}{3} \int_0^{+\infty} \frac{\pi(\mathrm{d}r)}{r}, \quad \mathbb{K}[Y_t] = \frac{M_4}{4} \int_0^{+\infty} \frac{\pi(\mathrm{d}r)}{r}, \tag{15}$$

and

$$\mathrm{AC}_Y(h) = \left( \int_0^{+\infty} \frac{\pi(\mathrm{d}r)}{r} \right)^{-1} \int_0^{+\infty} \frac{e^{-rh} \pi(\mathrm{d}r)}{r}. \tag{16}$$

These cumulants are effective statistics that determine the shape of probability density of discharge. They can be derived directly by taking suitable expectations of (1) as shown in the literature by virtue of its linearity (e.g., Barndorff-Nielsen and Stelzer, 2013; Yoshioka, 2022)[44,80]. These (unnormalized) cumulants are proportional to moments of $\nu$ and $\pi$ in a separable way, facilitating their identifications from data.

If $\pi$ has the gamma density $\gamma_{\alpha_r, \beta_r}$ with $\alpha_r > 1$ and $\beta_r > 0$, then the autocorrelation exhibits a sub-exponential decay:

$$\mathrm{AC}_Y(h) = \frac{1}{(1 + \beta_r h)^{\alpha_r - 1}}. \tag{17}$$

The process $X$ is said to exhibit truly long memory if $\mathrm{AC}_Y(h)$ is not integrable ($\alpha_r \in (1, 2]$), and moderately long memory if it is ($\alpha_r > 2$). In the sequel, we write $\overline{V} = \mathbb{V}[Y_t]$.

### 3.1.2 supOUSV process

We have the following proposition concerning the average and variance of $X$ and the covariance between $X$ and $Y$, which can be utilized in applications such as the moment matching for identifying a model.

***Proposition 1***

*For any $t \in \mathbb{R}$, the followings hold true:*

$$\mathbb{E}[X_t] = 0, \tag{18}$$

$$\mathbb{V}[X_t] = \underbrace{\frac{1}{2}\sigma^2 \overline{Y}}_{\text{Variance due to diffusion}}$$
$$+ \underbrace{\mu^2 \overline{V} \left( \int_0^{+\infty} \frac{\pi(\mathrm{d}r)}{r} \right)^{-1} \int_{R=0}^{R=+\infty} \int_{P=0}^{P=+\infty} \int_{r=0}^{r=+\infty} \frac{1}{r} \frac{PR}{P+R} \left( \frac{1}{P+r} + \frac{1}{R+r} \right) \pi(\mathrm{d}r) \rho(\mathrm{d}R) \rho(\mathrm{d}P)}_{\text{Variance due to discharge-dependent drift}}, \quad (19)$$

and

$$\mathbb{E}[X_t Y_t] = \mu \overline{V} \left( \int_0^{+\infty} \frac{\pi(\mathrm{d}r)}{r} \right)^{-1} \int_{R=0}^{R=+\infty} \int_{r=0}^{r=+\infty} \frac{R}{r(R+r)} \pi(\mathrm{d}r) \rho(\mathrm{d}R). \quad (20)$$

The first term on the right-hand side of (19) represents the variance due to diffusion, which is proportional to $\sigma^2$ and arises from stochastic volatility. By contrast, the second term in the right-hand side of (19) represents the variance due to the discharge-dependent drift, which is proportional to $\mu^2$ and independent of the sign of $\mu$. The integrals in (19)-(20) are numerically evaluated in this paper. According to (20), the sign of $\mu$ determines the sign of the covariance, providing a potential method for distinguishing the supOUSV model in real data because it is able to deal with both positive and negative correlations between discharge and WQI by properly specifying $\mu$.

We also obtain the closed-form expression for the autocorrelation $\mathrm{AC}_X(h)$ of $X$.

**Proposition 2**

*For any $t \in \mathbb{R}$ and $h \geq 0$, the autocorrelation $\mathrm{AC}_X(h)$ is expressed as follows:*

$$\mathrm{AC}_X(h) = \frac{1}{\mathbb{V}[X_t]} \left( \underbrace{\frac{\sigma^2 \overline{Y}}{2} I_1(h)}_{\text{Correlation due to diffusion}} + \underbrace{\mu^2 \overline{V} \left( \int_0^{+\infty} \frac{\pi(\mathrm{d}r)}{r} \right)^{-1} (I_2(h) + I_3(h))}_{\text{Correlation due to discharge-dependent drift}} \right), \quad (21)$$

*where*

$$I_1(h) = \int_{R=0}^{R=+\infty} e^{-Rh} \rho(\mathrm{d}R), \quad (22)$$

$$I_2(h) = \int_{R=0}^{R=+\infty} \int_{P=0}^{P=+\infty} \int_{r=0}^{r=+\infty} \frac{RP}{r(P+R)} \left( \frac{1}{P+r} + \frac{1}{R+r} \right) e^{-Ph} \pi(\mathrm{d}r) \rho(\mathrm{d}R) \rho(\mathrm{d}P), \quad (23)$$

*and*

$$I_3(h) = \int_{R=0}^{R=+\infty} \int_{P=0}^{P=+\infty} \int_{r=0}^{r=+\infty} \frac{RP}{r(R+r)(P-r)} \left( e^{-rh} - e^{-Ph} \right) \pi(\mathrm{d}r) \rho(\mathrm{d}R) \rho(\mathrm{d}P). \quad (24)$$

While **Proposition 2** formally represents $\mathrm{AC}_X(h)$, it does not address the existence or decay rate, which is analyzed below. In the rest of **Section 3.1.2**, we assume the following gamma distributions to quantify the tail behavior of $I_k$:

$$\pi(\mathrm{d}r) = \gamma_{\alpha_r,\beta_r}(r)\mathrm{d}r \text{ and } \rho(\mathrm{d}R) = \gamma_{\alpha_R,\beta_R}(R)\mathrm{d}R \text{ with } \alpha_r - 1, \alpha_R, \beta_r, \beta_R > 0. \quad (25)$$

We have the following estimate of $\mathrm{AC}_X(h)$ for large $h > 0$.

*Proposition 3*

*Assume (25) with*

$$(\alpha_r - 1)(\alpha_R + 2) > 1 \text{ if } \alpha_r \in (1, 2]. \quad (26)$$

*Then, it follows that* $\lim_{h \to +\infty} \mathrm{AC}_X(h) = 0$ *as* $h \to +\infty$. *Moreover, if* $\alpha_r > 2$, *then*

$$\mathrm{AC}_X(h) \leq \sigma^2 O(h^{-\alpha_R}) + \mu^2 \left\{ O(h^{-\alpha_R}) + O(h^{-(\alpha_r - 2)}) \right\}. \quad (27)$$

*If* $\alpha_r \in (1, 2]$, *then*

$$\mathrm{AC}_X(h) \leq \sigma^2 O(h^{-\alpha_R}) + \mu^2 \left\{ O(h^{-\alpha_R}) + O(h^{-\eta}) \right\} \quad (28)$$

*with*

$$\eta = \frac{(\alpha_r - 1)(\alpha_R + 2) - 1}{\alpha_R + 2} > 0. \quad (29)$$

According to **Propositions 2** and **3**, the term with the heaviest tail is $I_3$ when $\alpha_R$ is large at which the memory of river discharge (i.e., flood events) dominates. Specifically, the case $\mu \neq 0$ arises due to the coupling of the two superposed processes and is a unique feature of the proposed model, and the tail of $\mathrm{AC}_X(h)$ in this case possibly becomes the longer than that with $\mu = 0$. An engineering implication of **Proposition 3** is that autocorrelation may decay very slowly or even fail to converge to zero if the noise components of $X$ and $Y$ have sufficiently long memory, characterized by small $\alpha_r, \alpha_R$, as implied in (29).

*Remark 1* The relationships (19)–(20) are useful for determining the associated Orlicz spaces when quantifying the impact of model uncertainty (Yoshioka and Yoshioka, 2024a)[42] in the measures $\pi, \rho$ on statistics of $X$. A measure with stronger singularity (i.e., a more pronounced blow-up) at the origin corresponds to a smaller Orlicz space and is more sensitive to modeling errors.

### 3.2 Generalized Riccati equation

We have investigated the process $X$, while what directly contributes to water quality is its exponential, i.e., $e^X$. We therefore analyze when the exponential moment of the form $M(q) = \mathbb{E}\left[ e^{qX_t} \right]$ ($q \in \mathbb{R}$) exists. We focus on cases $q > 0$ because cases $q \leq 0$ are less non-trivial.

It is well known that, for a finite-dimensional system of affine SDEs, the moment-generating function is governed by a generalized Riccati equation as a nonlinear ordinary differential equation (e.g., Duffie et al., 2003)[45]. The following generalized Riccati equation is obtained from the finite-dimensional system (**Section A. 1 of Appendix**). In a stationary state, we have

$$M(q) = \exp\left(\lim_{t \to +\infty} \phi_t\right), \tag{30}$$

where

$$\frac{d\phi_t}{dt} = \int_0^{+\infty}\int_0^{+\infty}\left(e^{\psi_t(R)z} - 1\right)\nu(dz)\pi(dr) - \mu\int_0^{+\infty} P\omega_t(P)\rho(dP), \ t > 0, \tag{31}$$

$$\frac{\partial \omega_t(R)}{\partial t} = -R\omega_t(R), \ t, R > 0, \tag{32}$$

and

$$\frac{\partial \psi_t(r)}{\partial t} = -r\psi_t(r) + \mu q \int_0^{+\infty} \omega_t(P)P\rho(dP) + \frac{\sigma^2}{2}\int_0^{+\infty}\{\omega_t(P)\}^2 P\rho(dP), \ t, r > 0 \tag{33}$$

subject to initial conditions:

$$\omega_0(\cdot) = q \ \text{and} \ \psi_0(\cdot) = \phi_0 = 0. \tag{34}$$

By (33) and (34), we have

$$\psi_t(r) = \int_{P=0}^{P=+\infty}\int_{s=0}^{s=t} Pe^{-r(t-s)}\left(\mu q e^{-Pt} + \frac{\sigma^2 q^2}{2}e^{-2Pt}\right)ds\rho(dP). \tag{35}$$

To analyze the generalized Riccati equation in more detail, we assume the tempered stable model $\nu(dz) = a_1 e^{-a_2 z} z^{-(1+a_3)}dz$. To ensure $M(q)$ is well-defined, by (31) we need to have a bound $\psi_t(r) \leq a_2$ for all $t \geq 0$ and $r > 0$.

## Proposition 4

*Assume that* $\nu(dz) = a_1 e^{-a_2 z} z^{-(1+a_3)}dz$ *with* $a_1, a_2 > 0$ *and* $a_3 < 1$. *To ensure* $M(q) < +\infty$, *it is necessary to have*

$$\max\{\mu, 0\}q + \frac{\sigma^2 q^2}{4} \leq a_2 e. \tag{36}$$

We set $q_{\max}$ as the maximum $q$ satisfying (36), which is found as follows:

$$q_{\max} = \frac{2a_2 e}{\mu + \sqrt{\mu^2 + \sigma^2 a_2 e}}. \tag{37}$$

From (4), the moments $\mathbb{E}\left[C_t^k\right]$ of the concentration $C$ exist only if $k \leq q_{\max}$. Physically, higher-order moments are undefined for WQIs that fluctuate more (i.e., larger $\sigma$) and are more positively correlated with the discharge $Y$ (i.e., larger $\mu$).

The condition (36) is specific to the tempered stable $\nu$, and becomes unnecessary if we use the generalized model $\nu(dz) = a_1 e^{-a_2 z^{a_4}} z^{-(1+a_3)}dz$ with $a_1, a_2 > 0$, $a_3 < 1$, and $a_4 > 1$ (e.g., Grabchak, 2021)[81]. A drawback of the generalized tempered stable model is the increase in model parameters and the necessity for a more complicated numerical algorithm for simulating supOU processes. The case study in **Section 4** shows that the tempered stable model fails in applications when $q_{\max} < 1$; in such cases, the average $C$ does not exist. This issue can be efficiently avoided with a tempered stable $\nu$ by redefining the supOU process from (1) as follows:

$$Y_t = \int_{s=-\infty}^{s=t} \int_{z=0}^{z=+\infty} \int_{r=0}^{r=+\infty} z^{\frac{1}{1+\varepsilon}} e^{-r(t-s)} N(\mathrm{d}s, \mathrm{d}z, \mathrm{d}r), \ t \in \mathbb{R} \tag{38}$$

with a regularization parameter $\varepsilon > 0$. Accordingly, the generalized Riccati equation is modified as follows: (31) becomes

$$\frac{\mathrm{d}\phi_t}{\mathrm{d}t} = \int_0^{+\infty} \int_0^{+\infty} \left( \exp\left( \psi_t(r) z^{\frac{1}{1+\varepsilon}} \right) - 1 \right) \nu(\mathrm{d}z) \pi(\mathrm{d}r) - \mu \int_0^{+\infty} R \omega_t(R) \rho(\mathrm{d}R), \ t > 0, \tag{39}$$

while (32)–(33) remain unchanged. Moments and autocorrelation functions of the processes $X$ and $Y$ are also unchanged, except that $M_k$ ($k = 0,1,2,...$) must be formally replaced by $z^k$ by $z^{k/(1+\varepsilon)}$ in (2). An advantage of this regularization (38) is that it allows the use of simple numerical algorithms for simulating tempered stable subordinators (e.g., the rejection sampling method in Algorithm 0 of Kawai and Masuda (2011)[82]). Moreover, this regularization does not alter the overall structure of the proposed model but completely avoids the theoretical restriction in (36). Additionally, this regularization is qualitatively equivalent to employing a generalized tempered stable model in terms of moments, as shown by the following equality: for any $p > 0$,

$$\int_0^{+\infty} z^{\frac{p}{1+\varepsilon}} \frac{\exp(-a_2 z)}{z^{1+a_3}} \mathrm{d}z = (1+\varepsilon) \int_0^{+\infty} z^p \frac{\exp(-a_2 z^{1+\varepsilon})}{z^{1+(1+\varepsilon)a_3}} \mathrm{d}z, \tag{40}$$

where the right-hand side corresponds to the $p$ th moment of a generalized tempered stable model.

*Remark 2* Moment explosion phenomena have been reported using theoretical criteria for Markovian models (Deelstra et al., 2019)[83] and stochastic Volterra models (Gerhold et al., 2019)[84].

*Remark 3* We end this section with a remark of some difference between the present and existing models. The memory is often considered as a delay in dynamics system models (e.g., Arbi, 2022a; Arbi, 2022b)[85,86]. Our approach is seemingly different from these dynamical system models because the memory in our model means autocorrelation, while we expect that there is some linkage between our and these models because the delay is often modelled by some kernel function, which maybe decomposed into exponential functions in some cases. To the best of the authors' knowledge, there has been no research that applied the superposition approach to the coupled water quantity-quality dynamics. The proposed model reduces or related to a model without long memory (Calvani and Perona, 2023; Latella et al., 2024; Olson et al., 2021; Yoshioka and Yoshioka, 2023a; Yoshioka and Yoshioka, 2023b)[21-25] by choosing $\pi, \rho$ to be Dirac's deltas along with suitable modifications.

## 4. Application

### 4.1 Study site

The study site is the Kisuki point in the midstream of the Hii River (a Class-A River in Japan), located in the eastern part of Shimane Prefecture (**Figure 2**). The watershed of the whole Hii River system covers an area of 2,540 (km$^2$), with the main river stream extending 153 km. Approximately 500,000 people reside within the watershed[1]. According to Takeda (2023)[87], the catchment area of the Kisuki point is 451 km$^2$ and is predominantly forest (83.5%), with farmland (8.7%) and residential areas (0.7%). Theu, the Kisuki point has a rural catchment. The mainstream Hii River also contains two cascading brackish lakes, Lake Shinji and Lake Nakaumi. These lakes serve as Ramsar sites, providing unique habitats for migratory birds and functioning as important wetlands[2]. The biology (Kawaida et al., 2024)[88], hydrodynamics (Hafeez and Inoue, 2024)[89], and water quality (Kim et al., 2024)[90] of these brackish lakes have been studied. The water quantity and quality of the mainstream Hii River are critical in determining the characteristics of these lakes.

The water quality at Kisuki point has been investigated by Takeda (2023)[87], focusing on multiple indices, including total nitrogen (TN: August 20, 1991, to December 27, 2022), total phosphorus (TP: August 20, 1991, to December 27, 2022), and total organic carbon (TOC: April 5, 2005, to December 27, 2022). These indices were sampled and measured almost weekly from 1991 to 2022. Part of this dataset has been investigated by Yoshioka and Yoshioka (2024a)[42] using a superposition of diffusion processes without incorporating discharge data. This paper extends their analysis by including enhanced data from the 2022 survey, courtesy of Dr. Ikuo Takeda, the author of Takeda (2023)[87]. TN and TP are WQIs that determine environmental criteria for water use and eutrophication in aquatic systems (Atazadeh et al., 2024; Tian et al., 2024; Yuan and Paul, 2024)[91-93]. TOC is an effective trace element for carbon cycling, both locally in river environments and globally (Biedunkova and Kuznietsov, 2024; Starr et al., 2023)[94,95].

Another dataset used in this study is the weekly to monthly time series of dissolved silica (DSi, amount of Silicon in water samples) measured in a 1-km river reach, including Kisuki point, from May 22, 2018, to December 1, 2022. Part of this dataset has been investigated in previous studies (Yoshioka and Yoshioka, 2023a; Yoshioka and Yoshioka, 2023b)[24,25]. The data are treated as representative of the Kisuki point, assuming longitudinal homogeneity of water quality over this short reach. DSi is a conservative tracer with small correlation to discharge and has been used for tracking hydrological pathways (Das et al., 2022; Hachgenei et al., 2024; Jewell et al., 2023)[96-98]. For DSi, the data from May 2, 2022, to December 1, 2022, are excluded from autocorrelation analysis due to a 150-day gap from the just before sampling date.

The water samples were directly obtained from the Hii River. After that, the WQIs are measured using the following methods: TN is determined by UV spectrophotometry after potassium peroxodisulfate decomposition; TP by molybdenum blue spectrophotometry after potassium peroxodisulfate decomposition; TOC using a TOC analyzer (TOC-Vcsn, Shimadzu); and DSi by ICP emission spectrometry (ICPE9000, Shimadzu). Hourly river discharge data for Kisuki point[3] from January 1, 2017, at 01:00 to December 31, 2022, at 24:00 are also used in this study. Time-series data for discharge and WQIs are presented in **Figures 3-4**. The average, variance, and skewness of the discharge and each WQI are summarized in **Table 1**. The discharge is positively skewed due to the balance between persistent baseflow and flood events. TN, TP, and TOC are also positively skewed, while DSi is negatively skewed, though the magnitude of skewness for DSi is an order smaller for DSi than for the other three indices.

---

[1] MLIT, Japan. https://www.mlit.go.jp/river/toukei_chousa/kasen/jiten/nihon_kawa/0713_hiikawa/0713_hiikawa_00.html. Last Accessed on November 24, 2024. In Japanese.
[2] MLIT, Japan. https://www.mlit.go.jp/tagengo-db/en/H30-00978.html Last accessed on November 24, 2024. In Japanese.

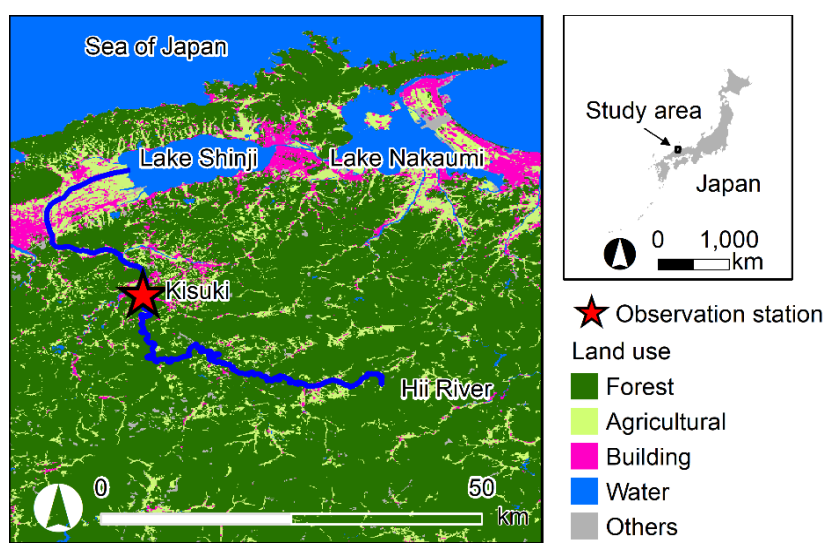

**Figure 2.** Study site and its surrounding area.

---

[3] MLIT, Japan. http://www1.river.go.jp/cgi-bin/SiteInfo.exe?ID=307041287705030 Last accessed on November 25, 2024. In Japanese.

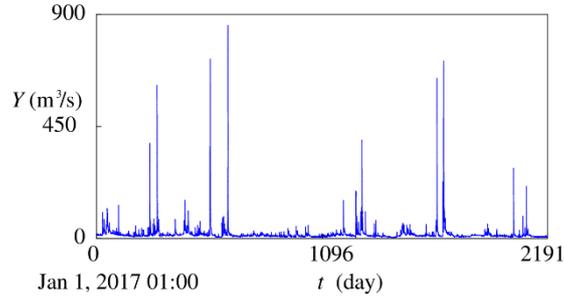

**Figure 3.** Hourly time-series data of discharge at Kisuki point.

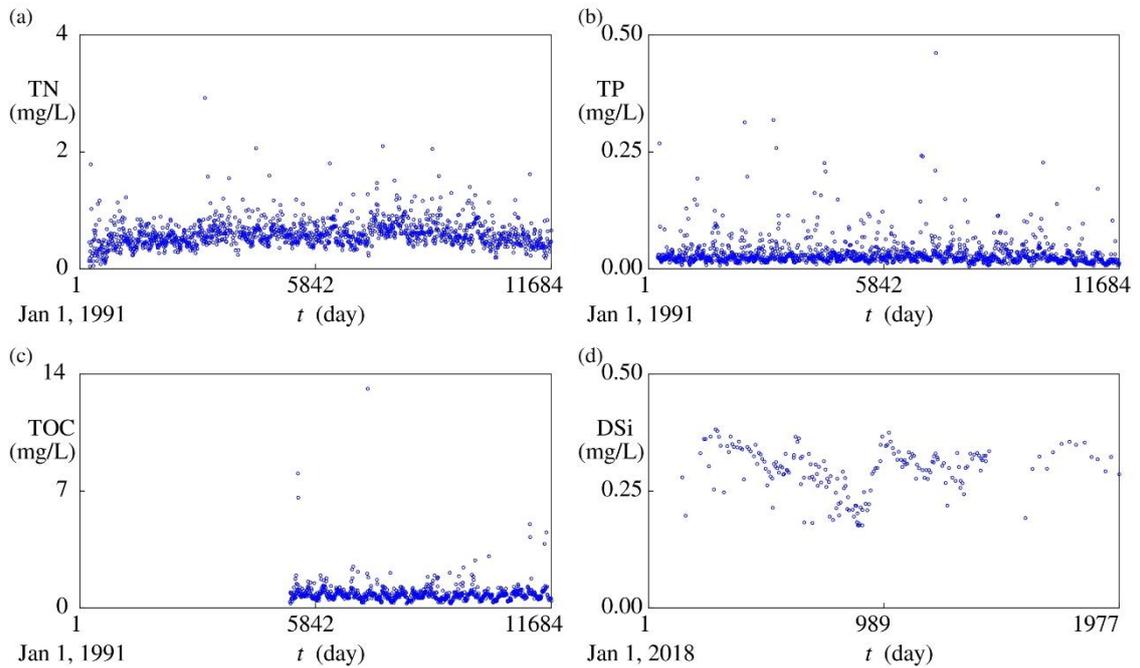

**Figure 4.** Time-series data of discharge and WQIs at Kisuki point: (a) TN, (b) TP, (c) TOC, and (d) DSi. Time shifts have been chosen for purely visualization purposes. TN=7.592 (mg/L) and TP=1.206 (mg/L) were recorded on 1993/6/29 (day 911), which are beyond the panels (a)-(b).

**Table 1.** Average, variance, and skewness of discharge and WQI. Units of variables are discharge ($m^3$/s), TN (mg/L), TP (mg/L), TOC (mg/L), and DSi (mg/L)

|          | Discharge   | TN         | TP         | TOC        | DSi         |
|----------|-------------|------------|------------|------------|-------------|
| Average  | 1.701.E+01  | 5.964.E-01 | 3.285.E-02 | 8.634.E-01 | 5.869.E+00  |
| Variance | 8.308.E+02  | 8.369.E-02 | 1.841.E-03 | 4.640.E-01 | 9.469.E-01  |
| Skewness | 1.406.E+01  | 1.056.E+01 | 1.574.E+01 | 1.014.E+01 | -7.141.E-01 |

## 4.2 Model identification

The supOUSV model for each WQI is identified as follows: We assume gamma distributions as in (25) and a tempered stable type with regularization, as in (38). The supOU process of discharge, which is common to all indices, is first identified using the moment-matching method (e.g., Yoshioka, 2022; Yoshioka and Yoshioka, 2024a)[42,80]. Here, the autocorrelation is fitted by minimizing the least-squares error between the empirical and theoretical values from (17) to identify $\pi$, and the Lévy measure $\nu$ is identified by minimizing the sum of relative errors in the average, variance, and skewness:

$$\left(\frac{\text{Ave}_e(Y) - \text{Ave}_m(Y)}{\text{Ave}_e(Y)}\right)^2 + \left(\frac{\text{Var}_e(Y) - \text{Var}_m(Y)}{\text{Var}_e(Y)}\right)^2 + \left(\frac{\text{Skew}_e(Y) - \text{Skew}_m(Y)}{\text{Skew}_e(Y)}\right)^2, \quad (41)$$

where Ave, Var, and Skew represent average, variance, and skewness of discharge $Y$, respectively. Subscripts "e" and "m" in (41) represent empirical and modeled values. Similar notations apply to $X$. The time lags for identifying autocorrelation are set at 30 days for discharge and 730 days for WQIs. We find that the identified models are not critically different when different lag durations are used (**Section A. 3**). We assume that the theoretical average of $X$ equals zero for each concentration (consistent with **Table 5** presented later).

After identifying the supOU process, the supOUSV model for each WQI, that is, $X$ is identified. The seasonality $\bar{C}e^{S_t}$ is determined by least-squares fitting between $\ln(\bar{C}e^{S_t})$ and observed concentrations $\ln C_t$, ensuring that the residual component $X_t$ has a numerically close-to-zero mean. Empirically, $n = 2$ is fixed to resolve the yearly variation of each WQI. Once seasonality component is identified, $e^{X_t}$ is calculated as $C_t / (\bar{C}e^{S_t})$, and hence $X_t$ by taking the logarithm. We need to identify the three quantities, which are the probability measure $\rho$ and parameters $\sigma$ and $\mu$. For DSi, which has a small correlation between $X$ and $Y$, we set $\mu = 0$ and fitted $\sigma$ using (19) as follows: $\sigma = \sqrt{\frac{2\text{Var}_e(X)}{\text{Ave}_e(Y)}}$. The measure $\rho$ is then fitted using a least-squares method based on (21) because, in this case, we have $\text{AC}_X(h) = (1 + \beta_R h)^{-\alpha_R}$ ($h \geq 0$).

For other indices (e.g., TN, TP, TOC), where the correlation between $X$ and $Y$ is not negligible $O(10^{-1})$, the above strategy does not apply. Instead, we used a two-step strategy. The temporal model assumption $\mu = 0$ is employed, and then we fit the auxiliary parameters $\alpha_{R,a}, \beta_{R,a}$ ($h \geq 0$). Next, using a trial-and-error method, we determine the parameter set $(\alpha_R, \sigma, \mu)$, as follows: Fix $\beta_R = \beta_{R,a}$ and rewrite the autocorrelation in (21) as

$$\text{AC}_X(h) = \frac{I_1(h) + w(I_2(h) + I_3(h))}{I_1(0) + w(I_2(0) + I_3(0))}, \quad h \geq 0 \quad (42)$$

with the weight

$$w = \frac{2\mu^2 \mathrm{Var}_e(Y)}{\sigma^2 \mathrm{Ave}_e(Y)} \left( \int_0^{+\infty} \frac{\pi(\mathrm{d}r)}{r} \right)^{-1}. \tag{43}$$

Here, we test several values of $\alpha_R$ around $\alpha_{R,a}$ (with an interval of 0.005 for TN and TP, and 0.05 for TOC). For each $\alpha_R$, we find $w$ such that the least-squares error between empirical and theoretical autocorrelation (42) is minimized. Having obtained the optimal $w$ for the given $\alpha_R$, we obtain $\sigma$ from

$$\mathrm{Var}_e(X) = \frac{1}{2} \sigma^2 \mathrm{Ave}_e(Y) \{ I_1(0) + w(I_2(0) + I_3(0)) \} \tag{44}$$

or equivalently

$$\sigma = \sqrt{\frac{2 \mathrm{Var}_e(X)}{\mathrm{Ave}_e(Y)\{I_1(0) + w(I_2(0) + I_3(0))\}}}, \tag{45}$$

Finally, we find $\mu$ from (43) as follows:

$$\mu = \frac{\mathbb{E}[X_t Y_t]}{|\mathbb{E}[X_t Y_t]|} \sigma \sqrt{\frac{w \mathrm{Ave}_e(Y)}{2 \mathrm{Var}_e(Y)} \int_0^{+\infty} \frac{\pi(\mathrm{d}r)}{r}} \tag{46}$$

ensuring that the signs of $\mu$ and $\mathbb{E}[X_t Y_t]$ are consistent (assuming that the latter is not zero). After obtaining $\sigma$ and $\mu$ for a given $\alpha_R$, the relative difference between empirical and theoretical correlations of $X$ and $Y$ is evaluated based on (20). The best parameter set is then selected by minimizing this difference. A key advantage of this method is that the theoretical variance exactly matches the empirical variance due to (45).

This empirical fitting strategy was chosen in this study due to direct application of continuous optimization methods for identifying $X$ is hindered by its complex and non-explicit dependence on $\alpha_R, \beta_R$. As such, the identified parameters in this paper are suboptimal. The integrals in the autocorrelation (21) are evaluated numerically using quantile-based quadrature (e.g., Yoshioka et al., 2023; Yoshioka and Yoshioka, 2024a)[42,70] and a forward-in-time explicit discretization is applied to the finite-dimensional version of the supOUSV model in **Section A. 1**. The computational resolution, which are specified by the degree-of freedoms $I_r, I_R$ in **Section A. 1**, are set as 2,048. The time increment of 0.01 day is used. Sampling $X$ of each WQI is conducted during a 1,000-year period after another 1,000-year period for burn-in with a deterministic initial condition.

### 4.3 Results and discussion

**Tables 2–3** summarize the identified parameter values. For discharge, models with ($\varepsilon = 0.1$) and without ($\varepsilon = 0$) regularization are examined. From **Table 2**, it is evident that the model without regularization ($\varepsilon = 0$) fails for all indices except DSi, due to $q_{\max} < 1$, which results in the nonexistence of the exponential moment of $X$ in the mean. This highlights the necessity of the proposed regularization method for this case study. Sensitivity of the discharge model to different values of $\varepsilon$ is explored in **Section A. 4** of **Appendix**, showing that the regularization introduced by $\varepsilon$ does not qualitatively affect the computational

results in this section. We therefore focus on the model with $\varepsilon = 0.1$. **Table 3** shows that the stochastic components of the four WQIs have long memories that would not be adequately captured by exponential ones. Particualrlly, those of TN and TP exhibit truly long memories.

In view of **Remark 1**, among the four WQIs, TN is the most vulnerable to modeling errors, followed by TP, TOC, and DSi. This is because smaller $\alpha_R$ implies less robustness. Moreover, evaluation of modeling uncertainty in TN should pay more attention than the other indices in this paper because the Orlicz space associated with TN is the smallest among the four indices in this paper (e.g., Proposition 1 of Yoshioka and Yoshioka (2024a))[42]. More specifically, assume that some modeling errors are included in the probability measure $\rho$ of the Gamam type as considered in this section. Then, an integrable function against $\rho$ should not have the singularity of the form $O(R^{-\gamma})$ with $\gamma \geq \alpha_R$. The theory of model uncertainty (Yoshioka and Yoshioka (2024a))[42] then suggests that the worst-case upper-bound of $X$ is more sensitive to the index with a smaller $\alpha_R$; namely, a process having a longer memory is more vulnerable to misspecification of this type.

**Table 4** compares the empirical and theoretical averages, variances, and skewnesses of discharge, showing strong agreement with relative errors smaller than 0.001%. **Table 5** shows the empirical statistics of $X$ for each WQI, where the regularized model is slightly more accurate for cumulants. **Figure 5** compares the empirical and theoretical autocorrelations of $Y$, computed from the identified models. Similarly, **Figure 6** compares the empirical and theoretical autocorrelations of $X$ (with $\mu \neq 0$ for TN, TP, and TOC). **Table 6** compares the empirical and theoretical covariances between $X$ and $Y$, suggesting that the identified model successfully captures the empirical correlations. Incorporating the discharge dependence in the SDE (11) through the parameter $\mu \neq 0$ is essential to obtain this result.

The empirical and computed probability density functions (PDFs) of $X$ for each WQI are compared in **Figure 7**. **Figure 8** shows the same comparison on a logarithmic scale to better visualize the PDF tails. **Figure 9** shows that the burn-in period of 1,000 years is sufficiently long so that fluctuations in the process $X$ is close to being stationary. As demonstrated in **Figures 7-8**, the identified models capture the empirical results, accurately representing the height and location of the PDF mode. The underestimation of the mode is due to the bounded range of empirical $X$ data, which is necessarily concentrated around the origin, while the model encompasses the entire range of $X$. Moreover, the model predicts the PDF tails, which are not always clear in the empirical data. The positive correlation between $X$ and discharge $Y$ manifests as heavier right tails in the theoretical PDFs for TN, TP, and TOC, whereas such behavior is absent for DSi.

Finally, we investigate the so-called concentration-discharge curve at flood events to reveal the complex hysteresis relationship between discharge and concentration of a WQI. We note that the water quality data used in this study were sampled weekly, and hence the direct use of the empirical data fails in this task. **Figure 10** shows some flood events and TN concentration generated by the Monte Carlo simulation and their concentration-discharge curve, which was extracted from a one-year sample path with the average of $C$ being 0.617 mg/L (empirical value is 0.596 mg/L) and the correlation between $Y$ and $C$ being 0.389 (empirical value is 0.441, and hence both are positive). In **Figure 10**, we plot the concentration-discharge curve by using daily-averaged $Y$ and $C$ so the shape of the curve becomes better visible. The resulting concentration-discharge curve is clockwise (**Figure 10(b)**) or counterclockwise (**Figure 10(d)**) depending on flood events. Event-based concentration-discharge curve directions have been reported for total reactive phosphorus around a Great Lakes (Dialameh and Ghane, 2023)[99], TN for a rural catchment (Rodríguez-Blanco et al., 2023)[100], dissolved organic carbon in a subtropical karst catchment (Qin et al., 2024)[58], and suspended sediment in the upper Yangtze River (Xue et al., 2024)[101]. Our model is consistent with these observations.

Physically, a concentration-discharge curve of a flood event in a river results from a complex hysteric dependence of the WQI not only on discharge but also on residual concentrations in surrounding groundwater of the river (Speir et al., 2024)[102]. Moreover, the correlation between water quantity and quality should also play a role; the WQIs considered in this paper have a positive correlation to the discharge, while some index, such as chloride and calcium ion, have a negative correlation (Cairoli et al., 2024; Knapp et al., 2022)[103,104]. Distance from the pollutant source to the river (Roberts et al., 2023)[105] and climatic and vegetation factors in watershed also affect concentration-discharge curves (Zhu et al., 2023)[106]. Such mechanisms are not directly considered in the proposed model; nevertheless, as demonstrated above, the proposed model could capture both the memory and correlation to discharge WQIs from the standpoint of infinite-dimensional SDEs, which is a novel contribution. Moreover, the proposed model, despite it is a simple model that does not account for such detailed physics, can also generate the two directions of the concentration-discharge curve found in the literature. To the best of the authors' knowledge, such a report does not exist for SDE models including long memory ones. Encoding a physics-informed event-based relationship into the proposed model requires the use of a higher-dimensional random measures in the SDEs. This would be a challenging and interesting topic to be resolved in future.

**Table 2.** Parameters for the model of discharge $Y$.

| | Not regularized ($\varepsilon = 0$) | Regularized ($\varepsilon = 0.1$) |
|---|---|---|
| $\alpha_r$ (-) | 2.143.E+00 | |
| $\beta_r$ (1/day) | 1.034.E+00 | |
| $a_1$ (m$^{3a_3+\frac{3\varepsilon}{1+\varepsilon}}$/s$^{a_3+\frac{\varepsilon}{1+\varepsilon}}$/day) | 1.266.E+00 | 1.124.E+00 |
| $a_2$ (s/m$^3$) | 1.960.E-03 | 8.920.E-04 |
| $a_3$ (-) | 8.084.E-01 | 7.500.E-01 |

**Table 3.** Parameter values for process $X$ of each WQI.

| | TN | TP | TOC | DSi |
|---|---|---|---|---|
| $\bar{C}$ (mg/L) | 5.553.E-01 | 2.642.E-02 | 7.676.E-01 | 5.745.E+00 |
| $A_1$ (-) | 7.104.E-02 | -2.752.E-01 | -2.512.E-01 | -1.133.E-01 |
| $A_2$ (-) | 6.562.E-02 | 7.639.E-02 | 6.511.E-02 | -2.330.E-02 |
| $B_1$ (-) | 7.198.E-01 | 1.344.E+00 | 5.396.E-01 | 1.801.E+00 |
| $B_2$ (-) | 7.185.E-01 | 8.146.E-01 | 1.640.E+00 | 4.360.E-01 |
| $\alpha_R$ (-) ($\mu = 0$) | 4.217.E-01 | 5.477.E-01 | 2.334.E+00 | 2.510.E+00 |
| $\beta_R$ (1/day) | 2.699.E-01 | 5.253.E-01 | 2.506.E-02 | 2.806.E-02 |
| $\alpha_R$ (-) ($\mu \neq 0$) | 3.750.E-01 | 4.850.E-01 | 2.650.E+00 | |
| $\sigma$ (s$^{1/2}$/m$^{3/2}$) | 1.077.E-01 | 1.483.E-01 | 1.037.E-01 | 5.483.E-02 |
| $\mu$ (s/m$^3$) | 2.752.E-02 | 2.917.E-02 | 2.567.E-02 | 0.000.E+00 |
| $q_{max}$ (-) | 1.899.E-01 | 1.768.E-01 | 2.033.E-01 | 2.663.E+00 |
| $w$ (1/day) | 7.532.E+00 | 4.469.E+00 | 7.072.E+00 | |

**Table 4.** Comparison between empirical and theoretical average, variance, and skewness of discharge $Y$.

| | Empirical | Theoretical ($\varepsilon = 0$) | Theoretical ($\varepsilon = 0.1$) |
|---|---|---|---|
| Average (m$^3$/s) | 1.701.E+01 | 1.701.E+01 | 1.701.E+01 |
| Variance (m$^6$/s$^2$) | 8.308.E+02 | 8.308.E+02 | 8.308.E+02 |
| Skewness (-) | 1.406.E+01 | 1.406.E+01 | 1.406.E+01 |

**Table 5.** Statistics of process $X$ for each WQI. Covariance and correlation between $X$ and $Y$ are calculated using daily-averaged $Y$.

| | TN | TP | TOC | DSi |
|---|---|---|---|---|
| Average of $X$ (-) | 5.890.E-07 | 3.634.E-06 | 2.628.E-07 | 1.120.E-06 |
| Variance of $X$ (-) | 1.373.E-01 | 2.794.E-01 | 1.459.E-01 | 2.556.E-02 |
| Skewness of $X$ (-) | -3.947.E-01 | 1.313.E+00 | 1.304.E+00 | -9.859.E-01 |
| Covariance between $X$ and $Y$ (m$^3$/s) | 2.671.E+00 | 5.161.E+00 | 2.436.E+00 | 8.105.E-02 |
| Correlation between $X$ and $Y$ (-) | 3.531.E-01 | 4.786.E-01 | 2.986.E-01 | 4.467.E-02 |

**Table 6.** Comparison between empirical and theoretical covariance (m$^3$/s) between $X$ (-) and $Y$ (m$^3$/s).

| | TN | TP | TOC | DSi |
|---|---|---|---|---|
| Empirical | 2.671.E+00 | 5.161.E+00 | 2.436.E+00 | 8.105.E-02 |
| Theoretical | 2.679.E+00 | 5.104.E+00 | 2.418.E+00 | 0.000.E+00 |

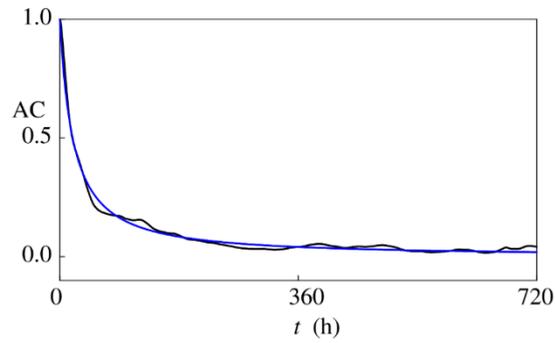

**Figure 5.** Comparison between empirical (black) and theoretical (blue) autocorrelation (AC) of discharge $Y$.

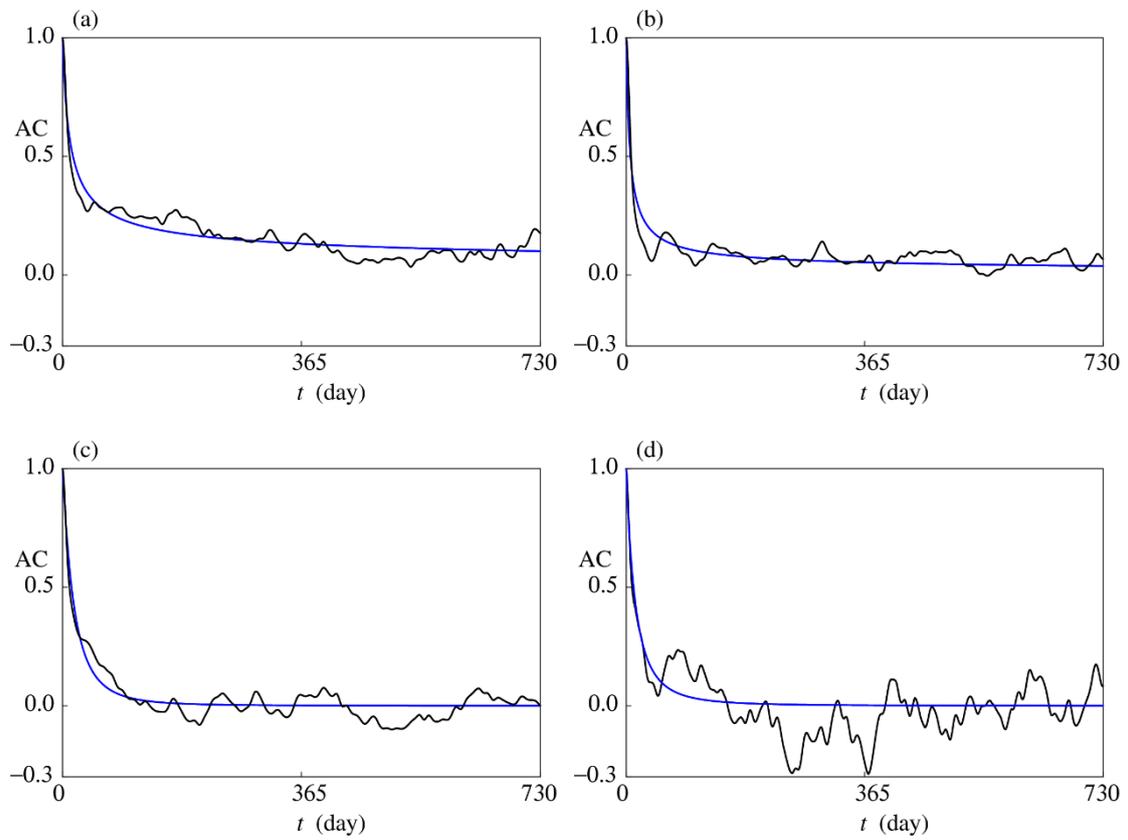

**Figure 6.** Comparison between empirical (black) and theoretical (blue) autocorrelations (AC) of $X$: (a) TN, (b) TP, (c) TOC, and (d) DSi.

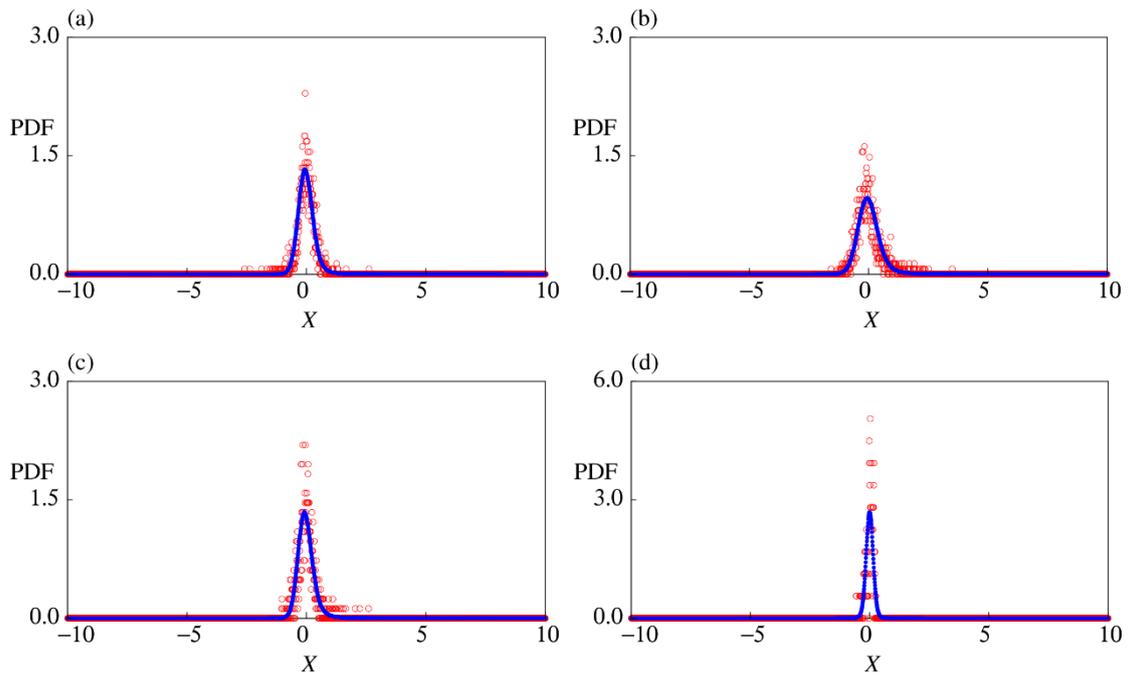

**Figure 7.** Comparison between empirical (red) and theoretical (blue) PDFs of $X$: (a) TN, (b) TP, (c) TOC, and (d) DSi.

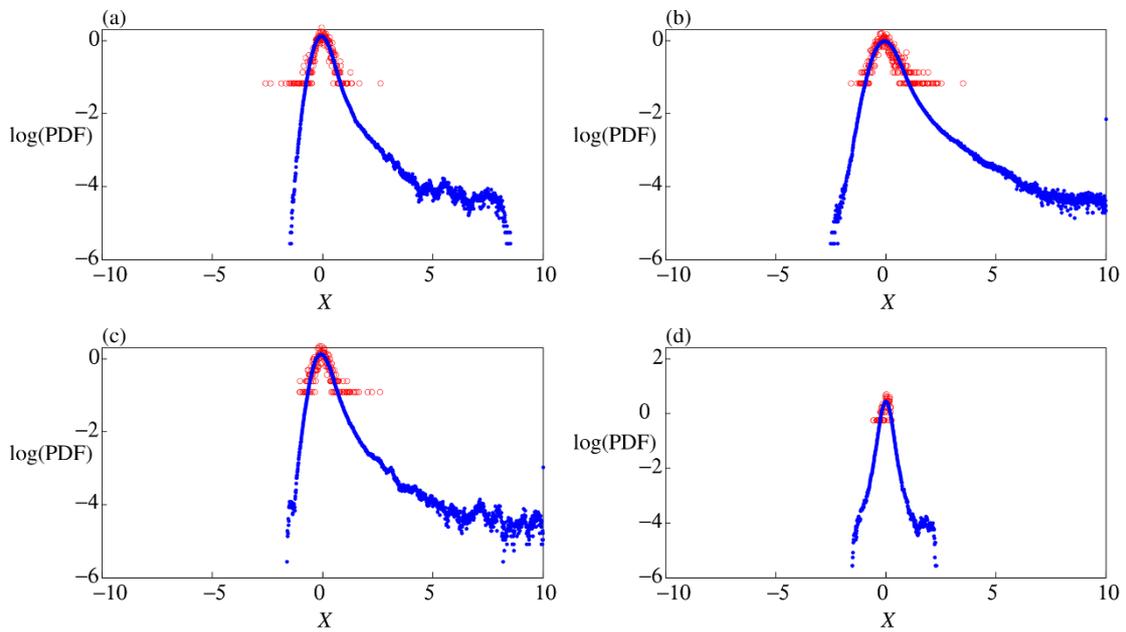

**Figure 8.** Comparison between empirical (red) and theoretical (blue) PDFs of $X$ in ordinary logarithmic scale: (a) TN, (b) TP, (c) TOC, and (d) DSi.

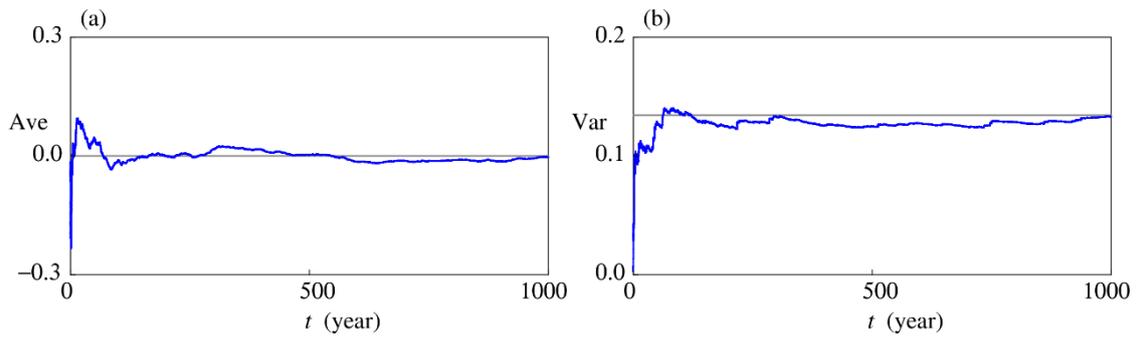

**Figure 9.** Convergence history of computed (a) average (Ave) and (b) variance (Var) of the process $X$ in the supOUSV model in the burn-in period. The grey line in each figure panel corresponds to theoretical value.

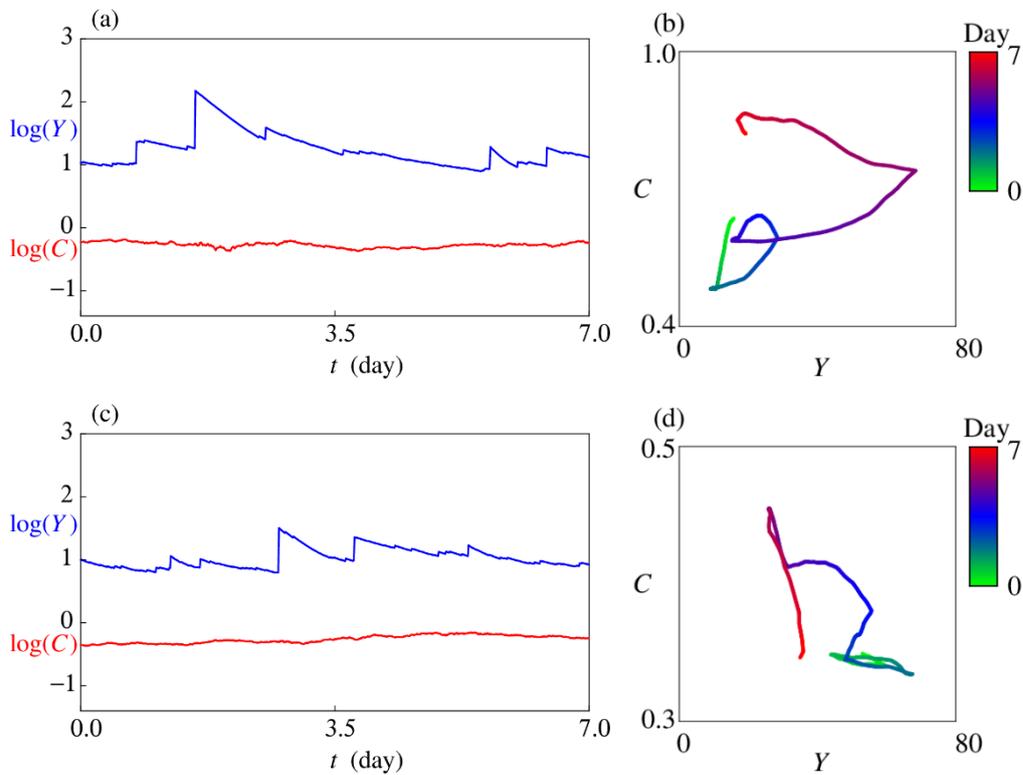

**Figure 10.** concentration-discharge relationships at flood events in model: (a) discharge $Y$ m$^3$/s and TN concentration $C$ mg/L during some simulated flood events in ordinary logarithmic scale and (b) corresponding concentration-discharge curve. Panels (c) and (d) show results for different flood events. Time is suitably normalized in panels (a) and (c), respectively.

## 5. Conclusion

We proposed a supOUSV model for coupled water quantity-quality dynamics in a river environment. The model demonstrated flexibility in representing the long-memory nature of discharge and that included in WQIs, as well as the ability to analytically obtain moments and autocorrelation, enabling the application of straightforward model identification methods. Analysis of the generalized Riccati equation provided the condition for defining each moment of the WQI for a specific case. The supOUSV model was applied to data collected at a study site in the midstream of the Hii River, clarifying its applicability and limitations for each WQI. The innovative integration of mathematical theory and computational techniques presented in this study would enhance the understanding of water environmental dynamics and their potential management strategies.

Several limitations of the supOUSV model are identified, which should be addressed in future research to develop a more accurate model with broader applicability. First, the effects of incorporating small processes with more generic forms have not been investigated and remain an open area of study. Another principal issue not addressed in this study is the biology and ecology of river environments, where water quantity and quality significantly impact aquatic organisms. The proposed model with the exponential nature of the seasonality decomposition allows for effectively separating the deterministic and stochastic parts of the positivity-preserving WQI dynamics, but its drawback is the possible loss of theoretical stability as indicated in **Proposition 3**. We will explore another decomposition method that does not lose the tractability of the model but is more stable. While space-time stochastic modeling could potentially address these interactions, the lack of spatially detailed water quality data near the study site made such modeling technically infeasible at this stage. Addressing this issue will require an intensive and large-scale survey to collect comprehensive data.

# Appendix

## A. 1 Finite-dimensional version

We present a finite-dimensional version of the supOUSV model and derive its associated generalized Riccati equation to enhance understanding of the infinite-dimensional versions discussed in the main text. The finite-dimensional version is formulated following Yoshioka et al. (2023)[75], who discussed a discretized supOU-type process.

Let $I_r, I_R \in \mathbb{N}$ be the degrees of freedom of the small processes of $X, Y$, respectively. Define nonnegative sequences $\{c_i\}_{1 \leq i \leq I_r}$ and $\{d_j\}_{1 \leq j \leq I_R}$ that weakly approximate the probability measures $\pi$ and $\rho$, respectively. We assume that $\sum_{i=1}^{I_r} c_i = \sum_{j=1}^{I_R} d_j = 1$. We prepare $I_r$ independent Poisson random measures $N^{(i)}$ ($1 \leq i \leq I_r$) with compensated measures $c_i \nu(\mathrm{d}z)\mathrm{d}s$. Also, we prepare $I_R$ independent Brownian motions $B^{(j)}$ ($1 \leq j \leq I_R$). All $N^{(i)}$ and $B^{(j)}$ are assumed to be mutually independent. Additionally, define discrete recession rates $\{r_i\}_{1 \leq i \leq I_r}$ and discrete reversion rates $\{R_j\}_{1 \leq j \leq I_R}$, both positive strictly increasing.

Under these conditions, the finite-dimensional supOU process $\hat{Y} = (\hat{Y}_t)_{t \in \mathbb{R}}$ is defined as follows:

$$\hat{Y}_t = \sum_{i=1}^{I_r} y_t^{(i)} \qquad (47)$$

with

$$y_t^{(i)} = \int_{s=-\infty}^{s=t} \int_{z=0}^{z=+\infty} z e^{-r_i(t-s)} N^{(i)}(\mathrm{d}s, \mathrm{d}z). \qquad (48)$$

Similarly, the discretized auxiliary process $\hat{X} = (\hat{X}_t)_{t \in \mathbb{R}}$ is defined as follows:

$$\hat{X}_t = \sum_{j=1}^{I_R} x_t^{(j)} \qquad (49)$$

with

$$x_t^{(j)} = \int_{s=-\infty}^{s=t} e^{-R_j(t-s)} \left\{ d_j R_j \mu (\hat{Y}_s - \bar{Y}) \mathrm{d}s + \sigma \sqrt{R_j d_j \hat{Y}_s} \mathrm{d}B_s^{(j)} \right\}. \qquad (50)$$

The corresponding SDEs for (48) and (50) are given by

$$\mathrm{d}y_t^{(i)} = -r_i y_t^{(i)} \mathrm{d}t + \int_{z=0}^{z=+\infty} z N^{(i)}(\mathrm{d}t, \mathrm{d}z) \qquad (51)$$

and

$$\mathrm{d}x_t^{(j)} = -R_j \left\{ x_t^{(j)} - d_j \mu (\hat{Y}_t - \bar{Y}) \right\} \mathrm{d}t + \sigma \sqrt{R_j d_j \hat{Y}_t} \mathrm{d}B_t^{(j)}, \qquad (52)$$

respectively.

According to Duffie et al. (2003)[45], SDEs (51) and (52) are affine. Therefore, it is possible to obtain the generalized Riccati equation associated with the moment-generating function $\hat{M}(q) = \mathbb{E}\left[e^{q\hat{X}_t}\right]$ in a stationary state. The generator $G$ associated with SDEs (51) and (52) is given by the linear partial integro-differential operator acting on a generic smooth function $\Phi = \Phi\left(x_1, x_2, ..., x_{I_R}, y_1, y_2, ..., y_{I_r}\right)$ with arguments in $\left(x_1, x_2, ..., x_{I_R}, y_1, y_2, ..., y_{I_r}\right) \in \mathbb{R}^{I_R + I_r}$ as follows:

$$G\Phi = \sum_{i=1}^{I_r}\left\{-r_i y_i \frac{\partial \Phi}{\partial y_i} + c_i \int_0^{+\infty}\left(\Phi\left(x_1, x_2, ..., x_{I_R}, y_1, y_2, ..., y_i + z_i, ...\right) - \Phi\right) \nu(dz_i)\right\}$$
$$+ \sum_{j=1}^{I_R}\left\{-R_j\left(x_j - \mu d_j\left(\sum_{i=1}^{I_r} y_i - \bar{Y}\right)\right)\frac{\partial \Phi}{\partial x_j} + \frac{1}{2}\sigma^2 R_j d_j\left(\sum_{i=1}^{I_r} y_i\right)\frac{\partial^2 \Phi}{\partial x_j^2}\right\} \quad (53)$$

Here, the first and second lines of (53) correspond to the generators of the discretized supOU and auxiliary processes, respectively. To estimate the functional form of the conditional moment-generating function (conditioned on the information at time 0)

$$\Phi\left(t, x_1, x_2, ..., x_{I_R}, y_1, y_2, ..., y_{I_r}\right) = \mathbb{E}\left[e^{q\hat{X}_t} \middle| \left(x_1, x_2, ..., x_{I_R}, y_1, y_2, ..., y_{I_r}\right)\right]$$
$$= \exp\left(\hat{\phi}_t + \sum_{i=1}^{I_r} \psi_t^{(i)} y_i + \sum_{j=1}^{I_R} \omega_t^{(j)} x_j\right), \quad (54)$$

we substitute (54) into the Kolmogorov equation:

$$\frac{\partial \Phi}{\partial t} = G\Phi, \ t > 0, \ \left(x_1, x_2, ..., x_{I_R}, y_1, y_2, ..., y_{I_r}\right) \in \mathbb{R}^{I_R + I_r}, \quad (55)$$

yielding (noting that $\Phi > 0$):

$$\frac{d}{dt}\left(\hat{\phi}_t + \sum_{i=1}^{I_r} \psi_t^{(i)} y_i + \sum_{j=1}^{I_R} \omega_t^{(j)} x_j\right)$$
$$= \sum_{i=1}^{I_r}\left\{-r_i y_i \psi_t^{(i)} + c_i \int_0^{+\infty}\left(e^{\psi_t^{(i)} z_i} - 1\right)\nu(dz_i)\right\} \quad (56)$$
$$+ \sum_{j=1}^{I_R}\left\{-R_j x_j \omega_t^{(j)} - \mu R_j d_j \bar{Y} \omega_t^{(j)} + \mu R_j d_j\left(\sum_{i=1}^{I_r} y_i\right)\omega_t^{(j)} + \frac{1}{2}\sigma^2 R_j d_j\left(\sum_{i=1}^{I_r} y_i\right)\left(\omega_t^{(j)}\right)^2\right\}$$

and obtain the finite-dimensional generalized Riccati equation:

$$\frac{d\hat{\phi}_t}{dt} = \sum_{i=1}^{I_r} c_i \int_0^{+\infty}\left(e^{\psi_t^{(i)} z} - 1\right)\nu(dz) - \mu \bar{Y} \sum_{j=1}^{I_R} d_j R_j \omega_t^{(j)}, \ t > 0, \quad (57)$$

$$\frac{d\omega_t^{(j)}}{dt} = -R_j \omega_t^{(j)}, \ 1 \leq j \leq I_R, \ t > 0, \quad (58)$$

and

$$\frac{d\psi_t^{(i)}}{dt} = -r_i \psi_t^{(i)} + \mu \sum_{j=1}^{I_R} d_j R_j \omega_t^{(j)} + \frac{1}{2}\sigma^2 \sum_{j=1}^{I_R} d_j R_j \left(\omega_t^{(j)}\right)^2, \ 1 \leq i \leq I_r, \ t > 0 \quad (59)$$

subject to initial conditions:

$$\omega_0^{(\cdot)} = q \text{ and } \psi_0^{(\cdot)} = \phi_0 = 0. \quad (60)$$

Finally, letting $I_r, I_R \to +\infty$ and $t \to +\infty$ formally derives (31)-(33).

## A. 2 Proofs

*Proof of Proposition 1*

For the average in (18), by (13) we formally have

$$
\begin{aligned}
\mathbb{E}[X_t] &= \mathbb{E}\left[\int_{R=0}^{R=+\infty}\int_{s=-\infty}^{s=t} e^{-R(t-s)}\left\{R\mu(Y_s-\bar{Y})\rho(\mathrm{d}R)\mathrm{d}s + \sigma\sqrt{RY_s}B(\mathrm{d}R,\mathrm{d}s)\right\}\right]\\
&= \mathbb{E}\left[\int_{r=0}^{R=+\infty}\int_{s=-\infty}^{s=t} e^{-R(t-s)}R\mu(Y_s-\bar{Y})\rho(\mathrm{d}R)\mathrm{d}s\right]\\
&\quad + \underbrace{\mathbb{E}\left[\int_{R=0}^{R=+\infty}\int_{s=-\infty}^{s=t} e^{-R(t-s)}\sigma\sqrt{RY_s}B(\mathrm{d}R,\mathrm{d}s)\right]}_{=0} \\
&= \int_{R=0}^{R=+\infty}\int_{s=-\infty}^{s=t} e^{-R(t-s)}R\mu\,\mathbb{E}[Y_s-\bar{Y}]\rho(\mathrm{d}R)\mathrm{d}s\\
&= \int_{R=0}^{R=+\infty}\int_{s=-\infty}^{s=t} e^{-R(t-s)}R\mu(\bar{Y}-\bar{Y})\rho(\mathrm{d}R)\mathrm{d}s\\
&= 0
\end{aligned}
\tag{61}
$$

because $Y$ is stationary and independent of $B$.

For variance from (19), because $B$ is a space-time Gaussian random measure and (18), and $Y$ is a nonnegative stationary process, we have

$$
\begin{aligned}
\mathbb{V}[X_t] &= \mathbb{E}[X_t^2]\\
&= \mathbb{E}\left[\left(\int_{R=0}^{R=+\infty}\int_{s=-\infty}^{s=t} e^{-R(t-s)}\left\{R\mu(Y_s-\bar{Y})\rho(\mathrm{d}R)\mathrm{d}s + \sigma\sqrt{RY_s}B(\mathrm{d}R,\mathrm{d}s)\right\}\right)^2\right]\\
&= \mathbb{E}\left[\begin{array}{l}\int_{R=0}^{R=+\infty}\int_{s=-\infty}^{s=t} e^{-R(t-s)}\left\{R\mu(Y_s-\bar{Y})\rho(\mathrm{d}R)\mathrm{d}s + \sigma\sqrt{RY_s}B(\mathrm{d}R,\mathrm{d}s)\right\}\\ \times\int_{P=0}^{P=+\infty}\int_{u=-\infty}^{u=t} e^{-P(t-u)}\left\{P\mu(Y_u-\bar{Y})\rho(\mathrm{d}P)\mathrm{d}u + \sigma\sqrt{PY_u}B(\mathrm{d}P,\mathrm{d}u)\right\}\end{array}\right]\\
&= \sigma^2\mathbb{E}\left[\int_{R=0}^{R=+\infty}\int_{s=-\infty}^{s=t}\int_{P=0}^{P=+\infty}\int_{u=-\infty}^{u=t}\sqrt{RP}e^{-R(t-s)-P(t-u)}\sqrt{Y_sY_u}B(\mathrm{d}R,\mathrm{d}s)B(\mathrm{d}P,\mathrm{d}u)\right]\\
&\quad + \mu^2\mathbb{E}\left[\int_{R=0}^{R=+\infty}\int_{s=-\infty}^{s=t}\int_{P=0}^{P=+\infty}\int_{u=-\infty}^{u=t} RPe^{-R(t-s)-P(t-u)}(Y_s-\bar{Y})(Y_u-\bar{Y})\rho(\mathrm{d}R)\rho(\mathrm{d}P)\mathrm{d}s\mathrm{d}u\right].\\
&= \sigma^2\int_{R=0}^{R=+\infty}\int_{s=-\infty}^{s=t} e^{-2R(t-s)}\mathbb{E}[Y_s]R\rho(\mathrm{d}R)\mathrm{d}s\\
&\quad + \mu^2\int_{R=0}^{R=+\infty}\int_{s=-\infty}^{s=t}\int_{P=0}^{P=+\infty}\int_{u=-\infty}^{u=t} RPe^{-R(t-s)-P(t-u)}\mathbb{E}[(Y_s-\bar{Y})(Y_u-\bar{Y})]\rho(\mathrm{d}R)\rho(\mathrm{d}P)\mathrm{d}s\mathrm{d}u\\
&= \sigma^2\bar{Y}\int_{R=0}^{R=+\infty}\int_{s=-\infty}^{s=t} e^{-2R(t-s)}R\rho(\mathrm{d}R)\mathrm{d}s\\
&\quad + \mu^2\bar{V}\int_{R=0}^{R=+\infty}\int_{s=-\infty}^{s=t}\int_{P=0}^{P=+\infty}\int_{u=-\infty}^{u=t} RPe^{-R(t-s)-P(t-u)}\mathrm{AC}_Y(|s-u|)\rho(\mathrm{d}R)\rho(\mathrm{d}P)\mathrm{d}s\mathrm{d}u
\end{aligned}
\tag{62}
$$

We have

$$
\int_{R=0}^{R=+\infty}\int_{s=-\infty}^{s=t} e^{-2R(t-s)}R\rho(\mathrm{d}R)\mathrm{d}s = \int_{R=0}^{R=+\infty}\frac{1}{2R}R\rho(\mathrm{d}R) = \frac{1}{2}.
\tag{63}
$$

Due to

$$
\mathrm{AC}_Y(|s-u|) = \left(\int_0^{+\infty}\frac{\pi(\mathrm{d}r)}{r}\right)^{-1}\int_{r=0}^{r=+\infty}\frac{1}{r}e^{-r|s-u|}\pi(\mathrm{d}r),
\tag{64}
$$

we also have

$$\int_{R=0}^{R=+\infty}\int_{s=-\infty}^{s=t}\int_{P=0}^{P=+\infty}\int_{u=-\infty}^{u=t} RP e^{-R(t-s)-P(t-u)} \mathrm{AC}_Y(|s-u|)\rho(\mathrm{d}R)\rho(\mathrm{d}P)\mathrm{d}s\mathrm{d}u$$

$$= \left(\int_0^{+\infty}\frac{\pi(\mathrm{d}r)}{r}\right)^{-1}\int_{R=0}^{R=+\infty}\int_{s=-\infty}^{s=t}\int_{P=0}^{P=+\infty}\int_{u=-\infty}^{u=t}\int_{r=0}^{r=+\infty}\frac{RP}{r}e^{-R(t-s)-P(t-u)-r|s-u|}\pi(\mathrm{d}r)\rho(\mathrm{d}R)\rho(\mathrm{d}P)\mathrm{d}s\mathrm{d}u \quad . \quad (65)$$

$$= \left(\int_0^{+\infty}\frac{\pi(\mathrm{d}r)}{r}\right)^{-1}\int_{R=0}^{R=+\infty}\int_{P=0}^{P=+\infty}\int_{r=0}^{r=+\infty}\frac{RP}{r}\left(\int_{s=-\infty}^{s=t}\int_{u=-\infty}^{u=t}e^{-R(t-s)-P(t-u)-r|s-u|}\mathrm{d}s\mathrm{d}u\right)\pi(\mathrm{d}r)\rho(\mathrm{d}R)\rho(\mathrm{d}P)$$

The inner integral in the last line of (65) is rewritten as follows:

$$\int_{s=-\infty}^{s=t}\int_{u=-\infty}^{u=t}e^{-R(t-s)-P(t-u)-r|s-u|}\mathrm{d}s\mathrm{d}u = \int_{s=0}^{s=+\infty}\int_{u=0}^{u=+\infty}e^{-Rs-Pu-r|s-u|}\mathrm{d}s\mathrm{d}u$$

$$= \underbrace{\int_{s=0}^{s=+\infty}\int_{u=s}^{u=+\infty}e^{-Rs-Pu-r(u-s)}\mathrm{d}u\mathrm{d}s}_{u>s} + \underbrace{\int_{u=0}^{u=+\infty}\int_{s=u}^{s=+\infty}e^{-Rs-Pu-r(s-u)}\mathrm{d}s\mathrm{d}u}_{s>u}$$

$$= \int_{s=0}^{s=+\infty}e^{-Rs+rs}\int_{u=s}^{u=+\infty}e^{-(P+r)u}\mathrm{d}u\mathrm{d}s + \int_{u=0}^{u=+\infty}e^{-Pu+ru}\int_{s=u}^{s=+\infty}e^{-(R+r)s}\mathrm{d}s\mathrm{d}u \quad . \quad (66)$$

$$= \int_{s=0}^{s=+\infty}e^{-Rs+rs}\frac{e^{-(P+r)s}}{P+r}\mathrm{d}s + \int_{u=0}^{u=+\infty}e^{-Pu+ru}\frac{e^{-(R+r)u}}{R+r}\mathrm{d}u$$

$$= \frac{1}{P+R}\left(\frac{1}{P+r}+\frac{1}{R+r}\right)$$

By combining (62)–(66), we obtain (19).

For the covariance in (20), by (64), we have

$$\mathbb{E}[X_t Y_t] = \mathbb{E}\left[\int_{R=0}^{R=+\infty}\int_{s=-\infty}^{s=t}e^{-R(t-s)}\left\{R\mu(Y_s-\bar{Y})\rho(\mathrm{d}R)\mathrm{d}s + \sigma\sqrt{RY_s}B(\mathrm{d}R,\mathrm{d}s)\right\}Y_t\right]$$

$$= \mathbb{E}\left[\int_{R=0}^{R=+\infty}\int_{s=-\infty}^{s=t}e^{-R(t-s)}R\mu(Y_s-\bar{Y})Y_t\rho(\mathrm{d}R)\mathrm{d}s\right]$$

$$+ \underbrace{\mathbb{E}\left[\int_{R=0}^{R=+\infty}\int_{s=-\infty}^{s=t}e^{-R(t-s)}\sigma\sqrt{RY_s}Y_t B(\mathrm{d}R,\mathrm{d}s)\right]}_{=0}$$

$$= \mu\int_{R=0}^{R=+\infty}\int_{s=-\infty}^{s=t}e^{-R(t-s)}R\mathbb{E}\left[Y_t Y_s - Y_t\bar{Y}\right]\rho(\mathrm{d}R)\mathrm{d}s$$

$$= \mu\int_{R=0}^{R=+\infty}\int_{s=-\infty}^{s=t}e^{-R(t-s)}R\mathbb{E}\left[Y_t Y_s - \bar{Y}^2\right]\rho(\mathrm{d}R)\mathrm{d}s \quad (67)$$

$$= \mu\int_{R=0}^{R=+\infty}\int_{s=-\infty}^{s=t}e^{-R(t-s)}R\mathbb{E}\left[(Y_s-\bar{Y})(Y_t-\bar{Y})\right]\rho(\mathrm{d}R)\mathrm{d}s$$

$$= \mu\bar{V}\int_{R=0}^{R=+\infty}\int_{s=-\infty}^{s=t}e^{-R(t-s)}\mathrm{AC}_Y(t-s)R\rho(\mathrm{d}R)\mathrm{d}s$$

$$= \mu\bar{V}\left(\int_0^{+\infty}\frac{\pi(\mathrm{d}r)}{r}\right)^{-1}\int_{R=0}^{R=+\infty}\int_{r=0}^{r=+\infty}\int_{s=-\infty}^{s=t}e^{-R(t-s)-r(t-s)}\frac{R}{r}\rho(\mathrm{d}R)\mathrm{d}s\pi(\mathrm{d}r)$$

$$= \mu\bar{V}\left(\int_0^{+\infty}\frac{\pi(\mathrm{d}r)}{r}\right)^{-1}\int_{R=0}^{R=+\infty}\int_{r=0}^{r=+\infty}\frac{R}{r}\left(\int_{s=-\infty}^{s=t}e^{-R(t-s)-r(t-s)}\mathrm{d}s\right)\pi(\mathrm{d}r)\rho(\mathrm{d}R)$$

and

$$\int_{s=-\infty}^{s=t}e^{-R(t-s)-r(t-s)}\mathrm{d}s = \int_{s=-\infty}^{s=0}e^{(R+r)s}\mathrm{d}s = \frac{1}{R+r} \quad . \quad (68)$$

Substituting (68) into (66), we obtain the expression for covariance (20).

□

*Proof of Proposition 2*

By the isometry and zero-mean Gaussian property of $B$ and $\mathbb{E}[X_t] = 0$, we have

$$\mathbb{V}[X_t]\text{AC}_X(h)$$
$$= \mathbb{E}[X_t X_{t+h}]$$
$$= \mathbb{E}\left[\begin{array}{l}\int_{R=0}^{R=+\infty}\int_{s=-\infty}^{s=t} e^{-R(t-s)}\left\{R\mu(Y_s - \bar{Y})\rho(\mathrm{d}R)\mathrm{d}s + \sigma\sqrt{RY_s}B(\mathrm{d}R,\mathrm{d}s)\right\}\\ \times\int_{P=0}^{P=+\infty}\int_{u=-\infty}^{u=t+h} e^{-P(t+h-s)}\left\{P\mu(Y_u - \bar{Y})\rho(\mathrm{d}P)\mathrm{d}u + \sigma\sqrt{PY_u}B(\mathrm{d}P,\mathrm{d}u)\right\}\end{array}\right]$$
$$= \sigma^2\mathbb{E}\left[\int_{R=0}^{R=+\infty}\int_{s=-\infty}^{s=t}\int_{P=0}^{P=+\infty}\int_{u=-\infty}^{u=t+h} \sqrt{RP}e^{-R(t-s)-P(t+h-u)}\sqrt{Y_s Y_u}B(\mathrm{d}R,\mathrm{d}s)B(\mathrm{d}P,\mathrm{d}u)\right]$$
$$+\mu^2\mathbb{E}\left[\int_{R=0}^{R=+\infty}\int_{s=-\infty}^{s=t}\int_{P=0}^{P=+\infty}\int_{u=-\infty}^{u=t+h} RPe^{-R(t-s)-P(t+h-u)}(Y_s - \bar{Y})(Y_u - \bar{Y})\rho(\mathrm{d}R)\rho(\mathrm{d}P)\mathrm{d}s\mathrm{d}u\right], \quad (69)$$
$$= \sigma^2\bar{Y}\int_{R=0}^{R=+\infty}\int_{s=-\infty}^{s=t} e^{-2R(t-s)-Rh}R\rho(\mathrm{d}R)\mathrm{d}s$$
$$+\mu^2\int_{R=0}^{R=+\infty}\int_{s=-\infty}^{s=t}\int_{P=0}^{P=+\infty}\int_{u=-\infty}^{u=t+h} RPe^{-R(t-s)-P(t+h-u)}\bar{V}\text{AC}_Y(|s-u|)\rho(\mathrm{d}R)\rho(\mathrm{d}P)\mathrm{d}s\mathrm{d}u$$
$$= \frac{\sigma^2\bar{Y}}{2}\int_{R=0}^{R=+\infty} e^{-Rh}\rho(\mathrm{d}R)$$
$$+\mu^2\bar{V}\int_{R=0}^{R=+\infty}\int_{s=-\infty}^{s=t}\int_{P=0}^{P=+\infty}\int_{u=-\infty}^{u=t+h} RPe^{-R(t-s)-P(t+h-u)}\text{AC}_Y(|s-u|)\rho(\mathrm{d}R)\rho(\mathrm{d}P)\mathrm{d}s\mathrm{d}u$$

where we used (63).

For the last integral in (69), we proceed as follows:

$$\int_{R=0}^{R=+\infty}\int_{s=-\infty}^{s=t}\int_{P=0}^{P=+\infty}\int_{u=-\infty}^{u=t+h} RPe^{-R(t-s)-P(t+h-u)}\text{AC}_Y(|s-u|)\rho(\mathrm{d}R)\rho(\mathrm{d}P)\mathrm{d}s\mathrm{d}u$$
$$= \left(\int_0^{+\infty}\frac{\pi(\mathrm{d}r)}{r}\right)^{-1}\int_{R=0}^{R=+\infty}\int_{s=-\infty}^{s=t}\int_{P=0}^{P=+\infty}\int_{u=-\infty}^{u=t+h}\int_{r=0}^{r=+\infty}\frac{RP}{r}e^{-R(t-s)-P(t+h-u)-r|s-u|}\pi(\mathrm{d}r)\rho(\mathrm{d}R)\rho(\mathrm{d}P)\mathrm{d}s\mathrm{d}u \quad (70)$$
$$= \left(\int_0^{+\infty}\frac{\pi(\mathrm{d}r)}{r}\right)^{-1}\int_{R=0}^{R=+\infty}\int_{P=0}^{P=+\infty}\int_{r=0}^{r=+\infty}\frac{RP}{r}e^{-Ph}\left(\int_{s=-\infty}^{s=t}\int_{u=-\infty}^{u=t+h} e^{-R(t-s)-P(t-u)-r|s-u|}\mathrm{d}s\mathrm{d}u\right)\pi(\mathrm{d}r)\rho(\mathrm{d}R)\rho(\mathrm{d}P)$$

and

$$\int_{s=-\infty}^{s=t}\int_{u=-\infty}^{u=t+h} e^{-R(t-s)-P(t-u)-r|s-u|}\mathrm{d}s\mathrm{d}u = \int_{s=-\infty}^{s=0}\int_{u=-\infty}^{u=h} e^{Rs+Pu-r|s-u|}\mathrm{d}s\mathrm{d}u$$
$$= \int_{s=0}^{s=+\infty}\int_{u=-h}^{u=+\infty} e^{-Rs-Pu-r|s-u|}\mathrm{d}s\mathrm{d}u \quad . \quad (71)$$
$$= \int_{s=0}^{s=+\infty}\int_{u=0}^{u=+\infty} e^{-Rs-Pu-r|s-u|}\mathrm{d}s\mathrm{d}u + \int_{s=0}^{s=+\infty}\int_{u=-h}^{u=0} e^{-Rs-Pu-r|s-u|}\mathrm{d}s\mathrm{d}u$$

The first term on the last line of (71) is identical to (66). For the second term, we focus on the case where $P \neq r$ in the integrand since the Lebesgue measure of the set $\{P = r\}$ is zero in $(0,+\infty)^2$:

$$\int_{s=0}^{s=+\infty}\int_{u=-h}^{u=0} e^{-Rs-Pu-r|s-u|}\mathrm{d}s\mathrm{d}u = \int_{s=0}^{s=+\infty}\int_{u=-h}^{u=0} e^{-Rs-Pu-r(s-u)}\mathrm{d}s\mathrm{d}u$$
$$= \int_{s=0}^{s=+\infty} e^{-(R+r)s}\mathrm{d}s\int_{u=-h}^{u=0} e^{-(P-r)u}\mathrm{d}u \quad . \quad (72)$$
$$= \frac{1}{R+r}\frac{1}{P-r}\left(e^{(P-r)h} - 1\right)$$

Substituting (71) and (72) into (70), we have

$$\int_{R=0}^{R=+\infty}\int_{s=-\infty}^{s=t}\int_{P=0}^{P=+\infty}\int_{u=-\infty}^{u=t+h} RPe^{-R(t-s)-P(t+h-u)}\mathrm{AC}_Y(|s-u|)\rho(\mathrm{d}R)\rho(\mathrm{d}P)\mathrm{d}s\mathrm{d}u$$

$$=\left(\int_0^{+\infty}\frac{\pi(\mathrm{d}r)}{r}\right)^{-1}\int_{R=0}^{R=+\infty}\int_{P=0}^{P=+\infty}\int_{r=0}^{r=+\infty}\frac{RP}{r}\frac{1}{P+R}\left(\frac{1}{P+r}+\frac{1}{R+r}\right)e^{-Ph}\pi(\mathrm{d}r)\rho(\mathrm{d}R)\rho(\mathrm{d}P). \quad (73)$$

$$+\left(\int_0^{+\infty}\frac{\pi(\mathrm{d}r)}{r}\right)^{-1}\int_{R=0}^{R=+\infty}\int_{P=0}^{P=+\infty}\int_{r=0}^{r=+\infty}\frac{RP}{r(R+r)(P-r)}\left(e^{-rh}-e^{-Ph}\right)\pi(\mathrm{d}r)\rho(\mathrm{d}R)\rho(\mathrm{d}P)$$

Finally, substituting (73) into (69) yields (21). Then, combining Equations (22), (23), and (24) directly follow from (21).

□

*Proof of Proposition 3*

Under the assumption of **Proposition 3**, we have

$$I_1(h)=\frac{1}{(1+\beta_R h)^{\alpha_R}}=O(h^{-\alpha_R}), \quad (74)$$

$$\begin{aligned} I_2(h) &\leq \int_{R=0}^{R=+\infty}\int_{P=0}^{P=+\infty}\int_{r=0}^{r=+\infty}\frac{RP}{r(P+R)}\left(\frac{1}{P+r}+\frac{1}{R+r}\right)e^{-Ph}\pi(\mathrm{d}r)\rho(\mathrm{d}R)\rho(\mathrm{d}P)\\ &\leq \int_{R=0}^{R=+\infty}\int_{P=0}^{P=+\infty}\int_{r=0}^{r=+\infty}\frac{RP}{r(P+R)}\left(\frac{1}{P}+\frac{1}{R}\right)e^{-Ph}\pi(\mathrm{d}r)\rho(\mathrm{d}R)\rho(\mathrm{d}P)\\ &=\left(\int_{r=0}^{r=+\infty}\frac{\pi(\mathrm{d}r)}{r}\right)\frac{1}{(1+\beta_R h)^{\alpha_R}}\\ &=O(h^{-\alpha_R}) \end{aligned} \quad (75)$$

and

$$\begin{aligned} I_3(h) &= \int_{R=0}^{R=+\infty}\int_{P=0}^{P=+\infty}\int_{r=0}^{r=+\infty}\frac{RP}{r(R+r)(P-r)}\left(e^{-rh}-e^{-Ph}\right)\pi(\mathrm{d}r)\rho(\mathrm{d}R)\rho(\mathrm{d}P)\\ &\leq \int_{R=0}^{R=+\infty}\int_{P=0}^{P=+\infty}\int_{r=0}^{r=+\infty}\frac{P}{r(P-r)}\left(e^{-rh}-e^{-Ph}\right)\pi(\mathrm{d}r)\rho(\mathrm{d}R)\rho(\mathrm{d}P)\\ &=\int_{P=0}^{P=+\infty}\int_{r=0}^{r=+\infty}\frac{P}{r(P-r)}\left(e^{-rh}-e^{-Ph}\right)\pi(\mathrm{d}r)\rho(\mathrm{d}P) \end{aligned} \quad (76)$$

We note the following elementary inequality, which will be used in the proof: for any $P>r>0$,

$$\frac{\partial}{\partial P}\left(\frac{e^{-rh}-e^{-Ph}}{P-r}\right)=\frac{(P-r)h+1-e^{(P-r)h}}{(P-r)^2}e^{-Ph}\leq 0 \quad (77)$$

due to $e^x\geq 1+x$ ($x\in\mathbb{R}$), and the limit

$$\lim_{P\to r}\frac{e^{-rh}-e^{-Ph}}{P-r}=he^{-rh},\; r>0. \quad (78)$$

The most technical aspect of this proof lies in evaluating the last integral in (76). We divide the domain of integration $(0,+\infty)\times(0,+\infty)$ as follows: for a small $\varepsilon>0$ specified later, $(0,+\infty)^2=D_1\cup D_2\cup D_3$ with $D_1=\{(P,r)\in(0,+\infty)^2:|P-r|\leq\varepsilon\}$, $D_2=\{(P,r)\in(0,+\infty)^2:P>r+\varepsilon\}$,

and $D_3 = \{(P, r) \in (0, +\infty)^2 : r > P + \varepsilon\}$. Then, from (77), for large $h > 0$, we have

$$\int_{D_1} \frac{P}{r(P-r)} \left(e^{-rh} - e^{-Ph}\right) \pi(dr) \rho(dP) \leq \int_{D_1} \frac{P}{r\varepsilon} \left(e^{-rh} - e^{-(r+\varepsilon)h}\right) \pi(dr) \rho(dP)$$
$$\leq \frac{1 - e^{-\varepsilon h}}{\varepsilon} \int_{D_1} \frac{P}{r} e^{-rh} \pi(dr) \rho(dP) \quad (79)$$
$$\leq \frac{1 - e^{-\varepsilon h}}{\varepsilon} \int_{P=0}^{P=+\infty} P \rho(dP) \int_{r=0}^{r=+\infty} \frac{1}{r} e^{-rh} \pi(dr)$$
$$= \frac{1 - e^{-\varepsilon h}}{\varepsilon} O\left(h^{1-\alpha_r}\right)$$

and

$$\int_{D_2} \frac{P}{r(P-r)} \left(e^{-rh} - e^{-Ph}\right) \pi(dr) \rho(dP) \leq \int_{D_2} \frac{P}{r\varepsilon} \left(e^{-Ph} - e^{-(P+\varepsilon)h}\right) \pi(dr) \rho(dP)$$
$$\leq \frac{1 - e^{-\varepsilon h}}{\varepsilon} \int_{D_2} \frac{P}{r} e^{-Ph} \pi(dr) \rho(dP) \quad (80)$$
$$\leq \frac{1 - e^{-\varepsilon h}}{\varepsilon} \int_{P=0}^{P=+\infty} P e^{-Ph} \rho(dP) \int_{r=0}^{r=+\infty} \frac{1}{r} \pi(dr)$$
$$= \frac{1 - e^{-\varepsilon h}}{\varepsilon} O\left(h^{-1-\alpha_R}\right)$$

From (77)-(78), we also have

$$\int_{D_3} \frac{P}{r(P-r)} \left(e^{-rh} - e^{-Ph}\right) \pi(dr) \rho(dP) \leq \int_{D_3} \frac{P}{r} h \max\left\{e^{-rh}, e^{-Ph}\right\} \pi(dr) \rho(dP)$$
$$\leq h \int_{D_3} \frac{P}{r} \left(e^{-rh} + e^{-Ph}\right) \pi(dr) \rho(dP) \quad (81)$$
$$= h \int_{D_3} \frac{P}{r} e^{-rh} \pi(dr) \rho(dP) + h \int_{D_3} \frac{P}{r} e^{-Ph} \pi(dr) \rho(dP)$$

Next, we estimate the last two integrals in (81):

$$\int_{D_3} \frac{P}{r} e^{-rh} \pi(dr) \rho(dP) = \int_{r=0}^{r=+\infty} \int_{P=\max\{0, r-\varepsilon\}}^{P=r+\varepsilon} \frac{P}{r} e^{-rh} \pi(dr) \rho(dP)$$
$$= \int_{r=0}^{r=+\infty} \frac{1}{r} e^{-rh} \left(\int_{P=\max\{0, r-\varepsilon\}}^{P=r+\varepsilon} P \rho(dP)\right) \pi(dr) \quad (82)$$
$$= O\left(h^{1-\alpha_r}\right) \times O\left(\varepsilon^{\alpha_R+1}\right)$$

and

$$\int_{D_3} \frac{P}{r} e^{-Ph} \pi(dr) \rho(dP) = \int_{P=0}^{P=+\infty} \int_{r=\max\{0, P-\varepsilon\}}^{r=P+\varepsilon} \frac{P}{r} e^{-Ph} \pi(dr) \rho(dP)$$
$$= \int_{P=0}^{P=+\infty} P e^{-Ph} \left(\int_{r=\max\{0, P-\varepsilon\}}^{r=P+\varepsilon} \frac{1}{r} \pi(dr)\right) \rho(dP) . \quad (83)$$
$$= O\left(h^{-1-\alpha_R}\right) \times O\left(\varepsilon^{\alpha_r-1}\right)$$

From (81)-(83), we obtain

$$\int_{D_3} \frac{P}{r(P-r)} \left(e^{-rh} - e^{-Ph}\right) \pi(dr) \rho(dP) \leq O\left(h^{2-\alpha_r} \varepsilon^{\alpha_R+1}\right) + O\left(h^{-\alpha_R} \varepsilon^{\alpha_r-1}\right). \quad (84)$$

Consequently, from (76), (79), (80), and (84), we arrive at the following estimate for large $h > 0$:

$$I_3(h) \leq \frac{1-e^{-\varepsilon h}}{\varepsilon} O\left(h^{-1-\alpha_R}\right) + \frac{1-e^{-\varepsilon h}}{\varepsilon} O\left(h^{1-\alpha_r}\right) + O\left(h^{2-\alpha_r}\varepsilon^{\alpha_R+1}\right) + O\left(h^{-\alpha_R}\varepsilon^{\alpha_r-1}\right). \tag{85}$$

Now, we choose $\varepsilon$ depending on $h$ as follows: $\varepsilon = h^{-\lambda}$ with some $\lambda \in (0,1)$. For large $h > 0$, we have

$$\frac{1-e^{-\varepsilon h}}{\varepsilon} = h^{\lambda}\left(1-\exp\left(-h^{1-\lambda}\right)\right) = O\left(h^{\lambda}\right), \tag{86}$$

and hence, since $\lambda \in (0,1)$,

$$\begin{aligned} I_3(h) &\leq O\left(h^{\lambda-1-\alpha_R}\right) + O\left(h^{\lambda+1-\alpha_r}\right) + O\left(h^{2-\alpha_r-\lambda(\alpha_R+1)}\right) + O\left(h^{-\alpha_R-\lambda(\alpha_r-1)}\right) \\ &\leq O\left(h^{-\alpha_R}\right) + O\left(h^{\lambda+1-\alpha_r}\right) + O\left(h^{2-\alpha_r-\lambda(\alpha_R+1)}\right) \end{aligned}. \tag{87}$$

Due to $\alpha_r > 1$, there exists some $\lambda = \lambda^* \in (0,1)$ such that $\lambda^* + 1 - \alpha_r < 0$. We consider two separate cases: $\alpha_r > 2$ and $\alpha_r \in (1,2]$. First, if $\alpha_r > 2$, we have, automatically,

$$2 - \alpha_r - \lambda^*(\alpha_R + 1) < 0. \tag{88}$$

By combining (74), (75), and (87), it follows that (without assuming (26)) $\lim_{h \to +\infty} \mathrm{AC}_X(h) = 0$ as $h \to +\infty$.

Second, if $\alpha_r \in (1,2]$, $\lambda^*$ such that $\lambda^* + 1 - \alpha_r < 0$ with (88) exists if

$$\frac{2-\alpha_r}{\alpha_R+1} < \alpha_r - 1 \Leftrightarrow 2 - \alpha_r < (\alpha_r-1)(\alpha_R+1) \Leftrightarrow 1 < (\alpha_r-1)(\alpha_R+2), \tag{89}$$

which is equivalent to the inequality in (26). From (26), we conclude that $\lim_{h \to +\infty} \mathrm{AC}_X(h) = 0$ as $h \to +\infty$ because we can choose $\lambda = \lambda^*$ to obtain

$$I_3(h) \leq O\left(h^{-\alpha_R}\right) + O\left(h^{\lambda^*+1-\alpha_r}\right) + O\left(h^{2-\alpha_r-\lambda^*(\alpha_R+1)}\right). \tag{90}$$

Finally, we close the proof by estimating $\mathrm{AC}_X(h)$. First, if $\alpha_r > 2$, we have, by taking $\lambda \to +0$,

$$I_3(h) \leq O\left(h^{-\alpha_R}\right) + O\left(h^{-(\alpha_r-2)}\right) \tag{91}$$

leading directly to (27). Second, if $\alpha_r \in (1,2]$, by equating the powers of $h$ in the second and third terms in (87), we have

$$I_3(h) \leq O\left(h^{-\alpha_R}\right) + O\left(h^{-\eta}\right) \tag{92}$$

with

$$\eta = \frac{(\alpha_r-1)(\alpha_R+2)-1}{\alpha_R+2} > 0 \tag{93}$$

due to Equation (26). This leads directly to the estimate in (29).

□

**Proof of Proposition 4**

We analyzed the upper bound of $\psi_t$. For this purpose, for each $r, P > 0$ with $P \neq r$, consider the following function of $t \geq 0$:

$$f_t = P\int_0^t e^{-r(t-s)} e^{-Ps} \, ds = \frac{P}{r-P}\left(e^{-Pt} - e^{-rt}\right). \tag{94}$$

We have $f_0 = 0$,

$$\frac{df_t}{dt} = \frac{P}{r-P}\left(-Pe^{-Pt} + re^{-rt}\right), \text{ and } \frac{d^2 f_t}{dt^2} = \frac{P}{r-P}\left(P^2 e^{-Pt} - r^2 e^{-rt}\right). \tag{95}$$

The function $f$ is maximized at $t = T$ such that the right-hand side of (95) vanishes. This $T$ is obtained as $T = \frac{1}{P-r}\ln\left(\frac{P}{r}\right)$. We then have the maximum value of the function $f$ as follows:

$$f_T = \frac{P}{r-P}\left(\left(\frac{P}{r}\right)^{-\frac{P}{P-r}} - \left(\frac{P}{r}\right)^{-\frac{r}{P-r}}\right) = \frac{P}{r-P}\left(\frac{P}{r}\right)^{-\frac{P}{P-r}}\left(1 - \left(\frac{P}{r}\right)^{\frac{P}{P-r}-\frac{r}{P-r}}\right) = \left(\frac{P}{r}\right)^{1-\frac{P/r}{P/r-1}}. \tag{96}$$

The right-hand side of (96) is expressed as $g(u)$ with $u = P/r > 0$, where $g(u) = u^{1-\frac{u}{u-1}} = u^{\frac{-1}{u-1}}$. We have $g(1) = e^{-1}$ by $g(u) = u^{1-\frac{u}{u-1}} = 1/(u-1+1)^{\frac{1}{u-1}}$ and taking the limit $u \to 1$. We have, for $u \neq 1$,

$$\frac{d}{du}\ln g(u) = \frac{d}{du}\left(\frac{1}{1-u}\ln u\right) = \frac{1}{(1-u)^2}\ln u - \frac{1}{u(1-u)} = \frac{1}{(1-u)^2}\left(\ln u + 1 - \frac{1}{u}\right). \tag{97}$$

We also have

$$\frac{d}{du}\left(\ln u + 1 - \frac{1}{u}\right) = \frac{1}{u} + \frac{1}{u^2} > 0 \text{ for } u > 0. \tag{98}$$

Moreover, it follows that

$$\lim_{u \to +0}\left(\ln u + 1 - \frac{1}{u}\right) = -\infty, \lim_{u \to +\infty}\left(\ln u + 1 - \frac{1}{u}\right) = +\infty, \text{ and } \ln 1 + 1 - \frac{1}{1} = 0. \tag{99}$$

This implies that $g(u)$ is maximized at $u = 1$ with the maximum value $g(1) = e^{-1}$.

The function $f$ with $P = r$ is given by $f_t = Pte^{-Pt}$, which is maximized at $t = P^{-1}$ as $e^{-1}$. Consequently, we have $0 \leq f_t \leq e^{-1}$ for any $P, r > 0$. We apply this to the right-hand side of (35): for all $t \geq 0$ and $r > 0$,

$$\begin{aligned}\psi_t(r) &= \int_{P=0}^{P=+\infty}\int_{s=0}^{s=t} Pe^{-r(t-s)}\left(\mu q e^{-Pt} + \frac{\sigma^2 q^2}{2}e^{-2Pt}\right)ds\rho(dP) \\ &= \mu q \int_{P=0}^{P=+\infty}\int_{s=0}^{s=t}\left(Pe^{-r(t-s)}e^{-Pt}\right)ds\rho(dP) + \frac{\sigma^2 q^2}{4}\int_{P=0}^{P=+\infty}\int_{s=0}^{s=t}\left(2Pe^{-r(t-s)}e^{-2Pt}\right)ds\rho(dP) \\ &\leq \max\{\mu,0\}q\int_{P=0}^{P=+\infty}\int_{s=0}^{s=t}\left(Pe^{-r(t-s)}e^{-Pt}\right)ds\rho(dP) + \frac{\sigma^2 q^2}{4}\int_{P=0}^{P=+\infty}\int_{s=0}^{s=t}\left(2Pe^{-r(t-s)}e^{-2Pt}\right)ds\rho(dP) \\ &\leq \left(\max\{\mu,0\}q + \frac{\sigma^2 q^2}{4}\right)e^{-1}\end{aligned} \tag{100}$$

Finally, from (100), we arrive at the strict upper bound $\psi_t(r) \leq a_2$ for all $t \geq 0$ and $R > 0$, provided condition in (36) is satisfied.

□

## A. 3 Identified values of $\alpha_r, \beta_r, \alpha_R, \beta_R$ for different time lag durations

For different time lag durations of autocorrelations, the identified values of the parameters $\alpha_r, \beta_r$ in the probability measure $\pi$ for the process $Y$ and $\alpha_r, \beta_r$ in the probability measure $\rho$ for the process $X$ are summarized in **Table A1**. The results demonstrate that the exponential decay of the autocorrelation in SS, the truly long-memory nature in TN and TP, and the moderately long-memory nature in TOC and DSi remain consistent, irrespective of the lag duration.

**Table A1.** Parameter values for process $X$ of each WQI. Numbers in "( )" represent lag durations used for parameter identification.

|  | TN | TP | TOC | DSi |
|---|---|---|---|---|
| $\alpha_R$ (-) (365 day) | 3.410.E-01 | 6.563.E-01 | 2.001.E+00 | 2.351.E+00 |
| $\beta_R$ (1/day) (365 day) | 4.970.E-01 | 3.376.E-01 | 3.105.E-02 | 3.056.E-02 |
| $\alpha_R$ (-) (730 day) | 4.217.E-01 | 5.477.E-01 | 2.334.E+00 | 3.510.E+00 |
| $\beta_R$ (1/day) (730day) | 2.699.E-01 | 5.253.E-01 | 2.506.E-02 | 2.806.E-02 |
| $\alpha_R$ (-) (1095 day) | 4.045.E-01 | 5.673.E-01 | 2.310.E+00 | 3.554.E+00 |
| $\beta_R$ (1/day) (1095 day) | 3.062.E-01 | 4.766.E-01 | 2.543.E-02 | 2.741.E-02 |

## A. 4 Sensitivity against $\varepsilon$

**Table A2** shows the identified parameter values and theoretical statistics of the proposed model for river discharge under different $\varepsilon$ values. This table serves as a supplementary version of **Table 2** in the main text. The results indicate that the average, variance, and skewness of the proposed model are not significantly affected by changes in $\varepsilon$, while kurtosis shows notable sensitivity. It is important to note that kurtosis is not included in the objective function in (41), and the autocorrelation of $X$ depends on neither skewness nor kurtosis, as established in **Proposition 2**. Moreover, the measure $\pi$, which governs the autocorrelation, remains independent of $\varepsilon$. Consequently, the regularization introduced by $\varepsilon$ does not qualitatively affect the computational results presented in the main text.

**Table A2.** Identified parameter values and theoretical statistics of model for river discharge against different values of $\varepsilon$.

|  | $\varepsilon = 0$ | $\varepsilon = 0.01$ | $\varepsilon = 0.1$ | $\varepsilon = 0.2$ | $\varepsilon = 0.4$ |
|---|---|---|---|---|---|
| $a_1$ ( $m^{3a_3+\frac{3\varepsilon}{1+\varepsilon}}$ / $s^{a_3+\frac{\varepsilon}{1+\varepsilon}}$ / day ) | 1.266.E+00 | 1.251.E+00 | 1.124.E+00 | 1.009.E+00 | 8.333.E-01 |
| $a_2$ (s/m$^3$) | 1.960.E-03 | 1.811.E-03 | 8.920.E-04 | 4.101.E-04 | 8.871.E-05 |
| $a_3$ (-) | 8.084.E-01 | 8.023.E-01 | 7.500.E-01 | 6.980.E-01 | 6.109.E-01 |
| Average (m$^3$/s) | 1.701.E+01 | 1.701.E+01 | 1.701.E+01 | 1.700.E+01 | 1.700.E+01 |
| Variance (m$^6$/s$^2$) | 8.308.E+02 | 8.308.E+02 | 8.308.E+02 | 8.308.E+02 | 8.308.E+02 |
| Skewness (-) | 1.406.E+01 | 1.406.E+01 | 1.406.E+01 | 1.406.E+01 | 1.406.E+01 |
| Kurtosis (-) | 4.089.E+02 | 4.074.E+02 | 3.948.E+02 | 3.831.E+02 | 3.652.E+02 |